\documentclass{amsart}
\usepackage{amsmath}
\usepackage{amssymb}

\newtheorem{theorem}{Theorem}
\theoremstyle{plain}

\newtheorem{claim}[theorem]{Claim}

\newtheorem{corollary}[theorem]{Corollary}

\newtheorem{definition}[theorem]{Definition}

\newtheorem{lemma}[theorem]{Lemma}

\newtheorem{proposition}[theorem]{Proposition}
\newtheorem{remark}[theorem]{Remark}

\numberwithin{equation}{section}
\numberwithin{theorem}{section}
\input{tcilatex}

\begin{document}
\title[Bounds via Self-Reducibility of the Weil Representation]{Bounds on
certain Higher-Dimensional Exponential Sums via the Self-Reducibility of the
Weil Representation}
\author{Shamgar Gurevich}
\curraddr{Department of Mathematics, University of California, Berkeley, CA
94720, USA. E-mail: shamgar@math.berkeley.edu}
\author{Ronny Hadani}
\curraddr{Department of Mathematics, University of Chicago, IL 60637, USA.
E-mail: hadani@math.uchicago.edu}

\begin{abstract}
We describe a new method to bound certain higher-dimensional exponential
sums which are associated with tori in symplectic groups over finite fields.
Our method is based on the self-reducibility property of the Weil
representation. As a result, we obtain a sharp form of the Hecke quantum
unique ergodicity theorem for generic linear symplectomorphisms of the $2N$%
-dimensional torus$.$
\end{abstract}

\maketitle

\section{Introduction}

\subsection{Eigenvectors of tori in the Weil representation}

To an $2N$-dimensional symplectic vector space ($V,\omega )$ over an odd
characteristic finite field $k=\mathbb{F}_{q}$ and a non-trivial additive
character $\psi :\mathbb{F}_{q}\rightarrow 
%TCIMACRO{\U{2102} }%
%BeginExpansion
\mathbb{C}
%EndExpansion
^{\ast }$ one can associate in a functorial manner a Hilbert space $\mathcal{%
H=H(}V)$ of dimension $q^{N}$ equipped with a unitary action of two symmetry
groups---the Heisenberg group $H=H(V)$ and the symplectic group $Sp=Sp(V).$
The first action is called the Heisenberg representation and we denote it by 
$\pi :H\rightarrow U(\mathcal{H)}$; the latter action is called the Weil
representation \cite{W} and we denote it by%
\begin{equation*}
\rho :Sp\rightarrow U(\mathcal{H)}.
\end{equation*}%
These representations play an important role in discrete harmonic analysis
with applications to various disciplines of pure and applied mathematics
such as representation theory, number theory, mathematical physics, coding
theory and signal processing.

Let $T\subset Sp$ be (the set of rational points of) a maximal torus which,
for simplicity, we will assume acts irreducibly on $V.$ The commutative
group $T$ acts, via the Weil representation, on the Hilbert space $\mathcal{H%
}$ and decomposes it into a direct sum of character spaces $\mathcal{%
H=\oplus H}_{\chi },$ where $\chi $ runs in the group of characters of $T.$
As it turns out $\dim (\mathcal{H}_{\chi })=1$ for each character $\chi $
that appears in the decomposition. In this paper we would like to study the
Wigner distribution 
\begin{equation}
\left \langle \varphi |\pi (v)\varphi \right \rangle ,  \label{Wigner}
\end{equation}%
associated with a unit character vector $\varphi \in \mathcal{H}_{\chi }$
and a non-zero vector $v\in V\subset H.$

\subsection{Wigner distributions as high-dimensional exponential sums}

We would like to bound the Wigner distribution. A possible method to achieve
this goal was developed in \cite{GH3, GH4}. The main idea is to write $%
\left
\langle \varphi |\pi (v)\varphi \right \rangle $ as an explicit
exponential sum for which we can use the powerful techniques of $\ell $-adic
cohomology to obtain good bounds.

If we denote by $P_{\chi }=\frac{1}{|T|}\sum_{g\in T}\chi ^{-1}(g)\rho (g)$
the orthogonal projector on the space $\mathcal{H}_{\chi }$ then $|T|\cdot
\left \langle \varphi |\pi (v)\varphi \right \rangle =|T|\cdot Tr(\pi
(v)P_{\chi })$ is equal to 
\begin{equation*}
c_{\chi }=\sum_{g\in T}\chi ^{-1}(g)\cdot Tr(\pi (v)\rho (g)),
\end{equation*}%
which by the formula \cite{GH6} for the trace of the Heisenberg--Weil
representation has the form of an explicit $N$-dimensional exponential sum
over $\mathbb{F}_{q}$%
\begin{equation*}
c_{\chi }=\sum_{g\in T\smallsetminus I}\chi ^{-1}(g)\cdot \sigma (\det
(g-I))\cdot \psi (\tfrac{1}{2}\omega (\frac{v}{g-I}{\small ,v})),
\end{equation*}%
where $\sigma :\mathbb{F}_{q}^{\ast }\rightarrow 
%TCIMACRO{\U{2102} }%
%BeginExpansion
\mathbb{C}
%EndExpansion
^{\ast }$ denote the Legendre character, $\psi :\mathbb{F}_{q}\rightarrow 
%TCIMACRO{\U{2102} }%
%BeginExpansion
\mathbb{C}
%EndExpansion
^{\ast }$ is a non-trivial additive character, and for the rest of this
paper\ $\frac{v}{g-I}=(g-I)^{-1}v.$

It is usually expected that these sums will have the \textquotedblleft
square root cancellation\textquotedblright \ phenomenon, i.e., that $%
|c_{\chi }|\leq d\sqrt{q}^{N}$ for some constant $d$ which is independent of 
$q.$ Using the formalism of $\ell $-adic cohomology, we replace the constant 
$c_{\chi }$ by the alternating sum 
\begin{equation}
c_{\chi }=\sum_{i=0}^{2N}(-1)^{i}Tr(Fr_{|H_{c}^{i}(\mathbf{X,}\mathcal{F}%
_{\chi })}),  \label{cohomolgy}
\end{equation}%
of the traces of the Frobenius operator acting on the cohomology groups with
compact support $H_{c}^{i}(\mathbf{X,}\mathcal{F}_{\chi })$ associated with
a suitable $\ell $-adic Weil sheaf $\mathcal{F}_{\chi }$ which lives on the
variety $\mathbf{X=T\smallsetminus I}$ where $\mathbf{T\subset Sp}$ is the
algebraic torus (we use boldface letters to denote algebraic varieties) such
that that $T=\mathbf{T(}\mathbb{F}_{q}).$

Deligne's Theory of weights \cite{D1} and a purity argument imply that the
eigenvalues of $Fr$ acting on $H_{c}^{i}(\mathbf{X,}\mathcal{F}_{\chi })$
are of absolute value $\sqrt{q}^{i}.$ Moreover, it can be shown that $%
H_{c}^{i}(\mathbf{X,}\mathcal{F}_{\chi })=0$ for $0\leq i\leq N-1.$ This
means that in order to obtain the expected bound on $c_{\chi }$ we need

\begin{itemize}
\item To show that all but the middle cohomology group vanish.

\item To calculate the dimension $d=\dim H_{c}^{N}(\mathbf{X,}\mathcal{F}%
_{\chi }).$
\end{itemize}

In \cite{GH3} we study the case $N=1$ and show that only the first
cohomology does not vanishes and $\dim H_{c}^{1}(\mathbf{X,}\mathcal{F}%
_{\chi })=2;$ therefore $|c_{\chi }|\leq 2\sqrt{q}.$ The computations for
general $N$ were carried in \cite{GH4}---indeed all but the middle
cohomology group vanish and $\dim H_{c}^{N}(\mathbf{X,}\mathcal{F}_{\chi
})=2^{N}.$ Hence, we find that in general 
\begin{equation}
|c_{\chi }|\leq 2^{N}\sqrt{q}^{N}.  \label{Es}
\end{equation}

However, as we will show in this paper (see the survey \cite{GH8} and the
announcement \cite{GH9}) the constant $2^{N}$ in the above bound is not
optimal---in fact we have%
\begin{equation}
|c_{\chi }|\leq 2\sqrt{q}^{N}.  \label{SE}
\end{equation}%
Thinking on the Frobenius operator acting on the space $H_{c}^{N}(\mathbf{X,}%
\mathcal{F}_{\chi })$ as a large matrix 
\begin{equation*}
Fr=%
\begin{pmatrix}
\lambda _{1} & \ast & \ast \\ 
& \ddots & \ast \\ 
&  & \lambda _{2^{N}}%
\end{pmatrix}%
\end{equation*}%
\smallskip a possible scenario which we might confront is cancellations
between different eigenvalues, more precisely angles, of the Frobenius
operator acting on a high-dimensional vector space, i.e., cancellations in
the sum $\sum_{j=1}^{2^{N}}e^{i\theta _{j}},$ where the angles $0\leq \theta
_{j}<2\pi $ are defined via $\lambda _{j}=e^{i\theta _{j}}\cdot \sqrt{q}%
^{N}. $ This problem is of a completely different nature, which is not
accounted for by standard cohomological techniques\footnote{%
We thank R. Heath-Brown for pointing out to us \cite{H} about the phenomenon
of cancelations between Frobenius eigenvalues in the presence of
high-dimensional cohomologies.}. \  \newline

\subsection{Sharp bound via self-reducibility}

Our approach to obtain the sharp bound (\ref{SE}) is to realize the constant 
$c_{\chi }$ as a one-dimensional exponential sum over $\mathbb{F}_{q^{N}}.$
This we do using a the self-reducibility property of the Weil representation.

\subsubsection{\textbf{Representation theoretic interpretation of the Wigner
distribution}}

The character vector $\varphi $ is a vector in a representation space $%
\mathcal{H}$ of the Weil representation of the symplectic group $Sp(2N,%
\mathbb{F}_{q}).$ The vector $\varphi $ is completely characterized in
representation theoretic terms, as being a character vector of the torus $T$%
. As a consequence, all quantities associated to $\varphi $, and in
particular the Wigner distribution $\left \langle \varphi |\pi (v)\varphi
\right \rangle $ is characterized in terms of the Weil representation. The
main observation to be made is that the character vector $\varphi $ can be
characterized in terms of another Weil representation, this time of a group
of a much smaller dimension---the torus $T$ induces an $\mathbb{F}_{q^{N}}$%
-structure on $(V,\omega )$ and now $\varphi $ is characterize in terms of
the Weil representation of $SL(2,\mathbb{F}_{q^{N}})$.

\subsubsection{\textbf{Self-reducibility property}}

A fundamental notion in our study is that of a \textit{symplectic module
structure}. A symplectic module structure is a triple $(K,V,\overline{\omega 
}),$ where $K$ is a finite dimensional commutative algebra over $k=\mathbb{F}%
_{q},$ equipped with an action on the vector space $V$, and $\overline{%
\omega }$ is a $K$-linear symplectic form satisfying the property $Tr_{K/k}(%
\overline{\omega })=\omega .$ Let $\overline{Sp}=Sp(V,\overline{\omega })$
be the group of $K$-linear symplectomorphisms with respect to the form $%
\overline{\omega }$. There exists a canonical embedding 
\begin{equation}
\iota :\overline{Sp}\hookrightarrow Sp.  \label{CI}
\end{equation}

It will be shown that associated to $T$ there exists a canonical symplectic
module structure $(K,V,\overline{\omega })$ so that $T\subset \overline{Sp}$%
. In our case the torus $T$ acts irreducibly on the vector space $V,$ hence,
the algebra $K$ is in fact a field with $\dim _{K}V=2$ which implies that $K=%
\mathbb{F}_{q^{N}}$ and $\overline{Sp}\simeq SL(2,\mathbb{F}_{q^{N}})$,
i.e., using (\ref{CI}) we get $T\subset SL(2,\mathbb{F}_{q^{N}})\subset Sp.$
Consider the Weil representation $(\rho ,Sp,\mathcal{H})$ associated with
the non-trivial additive character $\psi :\mathbb{F}_{q}\rightarrow 
%TCIMACRO{\U{2102} }%
%BeginExpansion
\mathbb{C}
%EndExpansion
^{\ast }$. Denote by $\overline{\psi }:K\rightarrow 
%TCIMACRO{\U{2102} }%
%BeginExpansion
\mathbb{C}
%EndExpansion
^{\ast }$ the additive character $\overline{\psi }=\psi \circ Tr_{K/k}.$

\begin{theorem}[\textbf{Self-reducibility property}]
\label{SRD}The restricted representation $(\overline{\rho }=\iota ^{\ast
}\rho ,\overline{Sp},\mathcal{H})$ is the Weil representation associated
with $\overline{\psi }$.
\end{theorem}

Applying the self-reducibility property to the torus $T$, it follows that
the vector $\varphi $ can be characterized in terms of the Weil
representation of $SL(2,\mathbb{F}_{q^{N}})$. Therefore, we can apply the
result obtained in \cite{GH3} and get the sharp bound $|c_{\chi }|\leq 2%
\sqrt{q}^{N}.$ Knowing that $|T|$ is of order of $q^{N}$, we obtain the
sharp bound on the Wigner distribution%
\begin{equation}
|\left \langle \varphi |\pi (v)\varphi \right \rangle |\leq \frac{2+o(1)}{%
\sqrt{q}^{N}},  \label{Sharp}
\end{equation}%
for every non zero vector $v\in V.$

We would like now to explain why the bounds on the Wigner distributions are
of interest.

\subsection{Quantum chaos problem}

One of the main motivational problems in quantum chaos is \cite{B1, B2, M, S}
describing eigenstates 
\begin{equation*}
\widetilde{H}\varphi =\lambda \varphi ,\allowbreak \  \  \varphi \in \mathcal{H%
},
\end{equation*}%
of a chaotic Hamiltonian $\widetilde{H}=Op(H):\mathcal{H\rightarrow H},$
where $\mathcal{H}$ is a Hilbert space. We deliberately use the notation $%
Op(H)$ to emphasize the fact that the quantum Hamiltonian $\widetilde{H}$ is
a quantization of a Hamiltonian $H:M\rightarrow 
%TCIMACRO{\U{2102} }%
%BeginExpansion
\mathbb{C}
%EndExpansion
$ where $M$ is a phase space---usually a cotangent bundle of a configuration
space $M=T^{\ast }X$, in which case $\mathcal{H=}L^{2}(X)$. In general,
describing $\varphi $ is considered to be an extremely complicated problem.
Nevertheless, for a few mathematical models of quantum mechanics rigorous
results have been obtained. We shall proceed to describe one of these models.

\subsubsection{\textbf{Hannay--Berry model}}

In \cite{HB} Hannay and Berry explored a model for quantum mechanics on the
two-dimensional symplectic torus $(\mathbb{T},\omega )$. Hannay and Berry
suggested to quantize simultaneously the functions on the torus and the
linear symplectic group $\Gamma \simeq SL(2,%
%TCIMACRO{\U{2124} }%
%BeginExpansion
\mathbb{Z}
%EndExpansion
)$. One of their main motivations was to study the phenomenon of quantum
chaos in this model \cite{M, R2}. More precisely, they considered an ergodic
discrete dynamical system on the torus which is generated by a hyperbolic
automorphism $A\in \Gamma $. Quantizing the system we replace---the phase
space $(\mathbb{T},\omega )$ by a finite dimensional Hilbert space $\mathcal{%
H}$; observables, i.e., functions $f\in C^{\infty }(\mathbb{T)}$ by
operators $\pi (f)\in End(\mathcal{H)}$; and symmetries by a unitary
representation $\rho :\Gamma \rightarrow U(\mathcal{H)}$ which, in
particular, enables one to associate to $A$ a unitary operator $\rho (A)$
acting on $\mathcal{H}$.

\subsubsection{The \textbf{Shnirelman theorem}}

Analogous with the case of the Schr\"{o}dinger equation, consider the
following eigenstates problem:

\begin{equation*}
\rho (A)\varphi =\lambda \varphi .
\end{equation*}

A fundamental result, valid for a wide class of quantum systems which are
associated to ergodic dynamics, is Shnirelman's theorem \cite{Sh}, asserting
that in the semi-classical limit \textquotedblleft almost
all\textquotedblright \ eigenstates become equidistributed in an appropriate
sense.

A variant of Shnirelman's theorem also holds in our situation \cite{BD}.
More precisely, we have that in the semi-classical limit $\hbar \rightarrow
0 $ for \textquotedblleft almost all\textquotedblright \ eigenstates $%
\varphi $ of the operator $\rho (A)$ the corresponding Wigner distribution $%
\left \langle \varphi |\pi (\cdot )\varphi \right \rangle :C^{\infty }(%
\mathbb{T})\rightarrow 
%TCIMACRO{\U{2102} }%
%BeginExpansion
\mathbb{C}
%EndExpansion
$ approaches the phase space average $\int_{\mathbb{T}}\cdot |\omega |$. In
this respect, it seems natural to ask whether there exist exceptional
sequences of eigenstates? Namely, eigenstates that do not obey the
Shnirelman's rule (\textquotedblleft scarred\textquotedblright \
eigenstates). It was predicted by Berry \cite{B1},\cite{B2}, that
\textquotedblleft scarring\textquotedblright \ phenomenon is not expected to
be seen for quantum systems associated with \textquotedblleft
generic\textquotedblright \ chaotic dynamics. However, in our situation the
operator $\rho (A)$ is not generic, and exceptional eigenstates were
constructed. Indeed, it was confirmed mathematically in \cite{FND} that
certain $\rho (A)$-eigenstates might localize. For example, in that paper a
sequence of eigenstates $\varphi $ was constructed, for which the
corresponding Wigner distribution approaches the measure $%
%TCIMACRO{\U{bd}}%
%BeginExpansion
{\frac12}%
%EndExpansion
\delta _{0}+%
%TCIMACRO{\U{bd}}%
%BeginExpansion
{\frac12}%
%EndExpansion
|\omega |$ on $\mathbb{T}$.

\subsubsection{\textbf{Hecke quantum unique ergodicity}}

A quantum system that obeys the Shnirelman rule is also called quantum
ergodic. Can one impose some natural conditions on the eigenstates so that
no exceptional eigenstates will appear? Namely, \textit{quantum unique
ergodicity} will hold. This question was addressed by Kurlberg and Rudnick
in \cite{KR1}, and they formulated a rigorous notion of Hecke quantum unique
ergodicity for the cases $\hbar =\frac{1}{p}$, $p$ a prime. Their basic
observation is that the degeneracies of the operator $\rho (A)$ are coupled
with the existence of symmetries---there exists a commutative group of
operators that commute with $\rho (A)$ and which can be computed
effectively. In more detail, the representation $\rho $ factors through the
Weil representation of the quotient group $Sp\simeq SL(2,\mathbb{F}_{p})$.
We denote by $T_{A}\subset Sp$ the centralizer of the element $A$, now
considered as an element of the quotient group. We call the group $T_{A}$
the \textit{Hecke torus} associated with $A$. The Hecke torus acts
semisimply on $\mathcal{H}$; therefore we have a decomposition into a direct
sum of Hecke eigenspaces $\mathcal{H=}\oplus \mathcal{H}_{\chi },$ where $%
\chi $ runs in the group of character of $T_{A}$. Consider a unit eigenstate 
$\varphi \in \mathcal{H}_{\chi }$ and the corresponding Wigner distribution $%
C^{\infty }(\mathbb{T})\mathbb{\rightarrow 
%TCIMACRO{\U{2102} }%
%BeginExpansion
\mathbb{C}
%EndExpansion
}$ defined by $f\mapsto \left \langle \varphi |\pi (f)\varphi \right \rangle
.$ The main statement in \cite{KR1} proves an explicit bound on the
semi-classical asymptotic---for sufficiently large $p$ they obtained $%
|\left
\langle \varphi |\pi (f)\varphi \right \rangle -\int_{\mathbb{T}%
}f|\omega ||\leq C_{f}/p^{1/4},$ where $C_{f}$ is a constant that depends
only on the function $f$. In addition, in \cite{R1, R2} Kurlberg and Rudnick
conjectured the following stronger bound:

\begin{equation}
\left \vert \left \langle \varphi |\pi (f)\varphi \right \rangle -\tint
\limits_{\mathbb{T}}f|\omega |\right \vert \leq \frac{C_{f}}{\sqrt{p}},
\label{estimate}
\end{equation}%
for sufficiently large prime $p.$

A particular case, which implies (\ref{estimate}), of the above inequality
is when $f=\xi $ a non-trivial character. In this case the integral $\int_{%
\mathbb{T}}\xi |\omega |$ vanishes and the bound (\ref{Sharp}) for the case
with $N=1$ and $k=\mathbb{F}_{p}$ gives $|\left \langle \varphi |\pi (\xi
)\varphi \right \rangle |\leq (2+o(1))/\sqrt{p}$, proving the conjecture 
\cite{GH3}.

\subsubsection{\textbf{The higher-dimensional Hannay--Berry model}}

The higher dimensional Hannay--Berry model is obtained as a quantization of
the $2N$-dimensional symplectic torus $(\mathbb{T},\mathbb{\omega )}$ acted
upon by the group $\Gamma \simeq Sp(2N,%
%TCIMACRO{\U{2124} }%
%BeginExpansion
\mathbb{Z}
%EndExpansion
)$ of linear symplectic autmorphisms. It was first constructed in \cite{GH2}%
, where, in particular, a quantization of the whole group of symmetries $%
\Gamma $ was obtained. Again, in the case $\hbar =\frac{1}{p}$ the
quantization of $\Gamma $ factors through the Weil representation of $%
Sp\simeq Sp(2N,\mathbb{F}_{p}).$ Considering a regular ergodic element $A\in
\Gamma $, i.e., $A$ \ generates an ergodic discrete dynamical system and it
is regular in the sense that it has distinct eigenvalues over $%
%TCIMACRO{\U{2102} }%
%BeginExpansion
\mathbb{C}
%EndExpansion
.$ It is natural to ask whether quantum unique ergodicity will hold true in
this setting as well, as long as one takes into account the whole group of
Hecke symmetries? Interestingly, the answer to this question is no. Several
new results in this direction have been announced recently. In the case
where the automorphism $A$ is \textit{non-generic}, meaning that it has an
invariant Lagrangian (and more generally co-isotropic) sub-torus $\mathbb{T}%
_{L}\subset \mathbb{T}$, an interesting new phenomenon was revealed. There
exists a sequence $\left \{ \varphi _{\hbar }\right \} $\textit{\ }of Hecke
eigenstates which might be related to the physical phenomenon of
\textquotedblleft localization\textquotedblright \ known in the literature
(cf. \cite{He},\cite{KH}) as \  \textquotedblleft scars\textquotedblright .
We will call them \textit{Hecke} scars. These states are localized in the
sense that the associated Wigner distribution converges to the Haar measure $%
\mu $\ on the invariant Lagrangian sub-torus 
\begin{equation}
\left \langle \varphi _{\hbar }|\pi (f)\varphi _{\hbar }\right \rangle
\rightarrow \tint \limits_{\mathbb{T}_{L}}f\mu ,\text{ as }\hbar \rightarrow
0,  \label{localization}
\end{equation}%
for every smooth observable $f$. These special Hecke eigenstates were first
established in \cite{Gu}. The semi-classical interpretation of the
localization phenomena (\ref{localization}) \ was announced in \cite{Ke}.

The above phenomenon motivates the following definition:

\begin{definition}
\label{qergodic}We will call an element $A\in \Gamma $ \textit{generic}%
\textbf{\ }if it is regular and admits no non-trivial invariant co-isotropic
sub-tori.
\end{definition}

\begin{remark}
The collection of generic elements constitutes an open subscheme of $\Gamma
. $ In particular, a generic element need not be ergodic automorphism of $%
\mathbb{T}$. However, in the case where $\Gamma \simeq SL_{2}(%
%TCIMACRO{\U{2124} }%
%BeginExpansion
\mathbb{Z}
%EndExpansion
)$ every ergodic (i.e., hyperbolic) element is generic. An example of
generic elements which is not ergodic is given by the Weyl element $w=%
{\small (}%
\begin{array}{cc}
0 & 1 \\ 
-1 & 0%
\end{array}%
{\small )}.$
\end{remark}

For the sake of simplicity let us assume now that the automorphism $A$ is 
\textit{strongly generic}, i.e., it has no non-trivial invariant sub-tori.
This case was first considered in \cite{GH4}, where using the bound (\ref{Es}%
) we obtain that for a fixed non-trivial character $\xi $ of $\mathbb{T}$ 
\begin{equation}
\left \vert \left \langle \varphi |\pi (\xi )\varphi \right \rangle \right
\vert \leq \frac{m_{\chi }\cdot (2+o(1))^{N}}{\sqrt{p}^{N}},
\label{NonSharp}
\end{equation}%
for a sufficiently large prime number $p$, where $m_{\chi }=\dim \mathcal{H}%
_{\chi }.$

In particular, using the bound (\ref{NonSharp}) we have:

\begin{theorem}[Hecke quantum unique ergodicity]
\label{QUE}Consider an observable $f$ $\in C^{\infty }(\mathbb{T)}$ and a
sufficiently large prime number $p.$ Then 
\begin{equation*}
\left \vert \left \langle \varphi |\pi (f)\varphi \right \rangle -\tint
\limits_{\mathbb{T}}fd\mu \right \vert \leq \frac{C_{f}}{\sqrt{p}^{N}},
\end{equation*}%
where $\mu =|\omega |^{N}$ is the corresponding volume form and $C_{f}$ is
an explicit computable constant which depends only on the function $f.$
\end{theorem}

The new method, using the self-reducibility property applied to the torus $%
T_{A},$ leads to a bound similar to (\ref{Sharp}) and to a sharp form of
Theorem\  \ref{QUE}.

\begin{theorem}[Sharp bound]
\textbf{\label{rate theorem}}Let $\xi $ be a non-trivial character of $%
\mathbb{T}$. For sufficiently large prime number $p$ the following bound
holds: 
\begin{equation}
\left \vert \left \langle \varphi |\pi (\xi )\varphi \right \rangle \right
\vert \leq \frac{m_{\chi }\cdot (2+o(1))^{r_{p}}}{\sqrt{p}^{N}},
\label{simplified3}
\end{equation}%
where the number $r_{p}$ is an integer between $1$ and $N$ that we call the
symplectic rank of $T_{A}$.
\end{theorem}

\begin{remark}
If the torus $T_{A}$ acts irreducibly on $V\simeq \mathbb{F}_{p}^{2N}$ then $%
r_{p}=1$; and if it splits, i.e., $T_{A}\simeq \mathbb{F}_{p}^{\ast N}$ then 
$r_{p}=N.$ In general $\left( \text{see Subsection \ref{St}}\right) $ the
distribution of the symplectic rank $r_{p\text{ }}$in the set $\left \{
1,...,N\right \} $ is governed by the Chebotarev density theorem applied to
a suitable Galois group $G$. For example, in the case where $A\in Sp(4,%
%TCIMACRO{\U{2124} }%
%BeginExpansion
\mathbb{Z}
%EndExpansion
)$ is strongly generic then $G$ is the symmetric group $S_{2}$ and we have
the density law \ 
\begin{equation*}
\lim_{x\rightarrow \infty }\frac{\# \{r_{p}=r\text{ }|\text{ }p\leq x\}}{\pi
(x)}=\tfrac{1}{2},\text{ \  \  \  \  \  \ }r=1,\text{ }2,
\end{equation*}%
where $\pi (x)$ denotes the number of primes up to $x$.
\end{remark}

\subsection{Quantum unique ergodicity for statistical states}

As in harmonic analysis, we would like to use Theorem \ref{rate theorem}
concerning the Hecke eigenstates in order to extract information on the
spectral theory of the operator $\rho (A)$ itself. For the sake of
simplicity, let us assume again that $A$ is strongly generic, i.e., it acts
on the torus $\mathbb{T}$ with no non-trivial invariant sub-tori. The
following is a possible reformulation of the quantum unique ergodicity
statement---one which is formulated for the automorphism $A$ itself instead
of the all Hecke group of symmetries. The element $A$ acts via the Weil
representation $\rho $ on the space $\mathcal{H}$ and decomposes it into a
direct sum of $\rho \left( A\right) $-eigenspaces 
\begin{equation}
\mathcal{H=}\tbigoplus \mathcal{H}_{\lambda }.  \label{dsA}
\end{equation}

Considering an $\rho \left( A\right) $-eigenstate $\varphi $ and the
corresponding projector $P_{\varphi }$ one usually studies the Wigner
distribution $\left \langle \varphi |\pi (f)\varphi \right \rangle =Tr(\pi
(f)P_{\varphi })$ which, due to the fact that $rank(P_{\varphi })=1,$ is
sometimes called$\ $a \textit{\textquotedblleft pure state\textquotedblright
.} In the same way, we might think about a Hecke--Wigner distribution $%
\left
\langle \varphi |\pi (f)\varphi \right \rangle =Tr(\pi (f)P_{\chi }),$
attached to a $T_{A}$-eigenstate $\varphi $, as a \textquotedblleft pure
Hecke state\textquotedblright . Following von Neumann \cite{vN} we suggest
the possibility of looking at the more general \textit{\textquotedblleft
statistical state\textquotedblright } defined by a non-negative self-adjoint
operator $D$---called the von Neumann density operator--- normalized to have 
$Tr(D)=1.$ For example, to the automorphism $A$ we can attach the natural
family of density operators $D_{\lambda }=\frac{1}{m_{\lambda }}P_{\lambda
,} $where $P_{\lambda }$ is the orthogonal projector on the eigenspace $%
\mathcal{H}_{\lambda }$ $\left( \text{\ref{dsA}}\right) ,$ and $m_{\lambda
}=\dim (\mathcal{H}_{\lambda }).$ Consequently, we obtain a family of
statistical states 
\begin{equation*}
Tr(\pi (\cdot )D_{\lambda }).
\end{equation*}

\begin{theorem}
\label{QUESS} Let $\xi $ be a non-trivial character of $\mathbb{T}$. For a
sufficiently large prime number $p$ we have 
\begin{equation}
|Tr(\pi (\xi )D_{\lambda })|\leq \frac{m\cdot (2+o(1))^{r_{p}}}{\sqrt{p}^{N}}%
,  \label{BQUESS}
\end{equation}%
where $1\leq r_{p}\leq N$ is an integer which is determined by $A,$ and $%
m=\max \dim \mathcal{H}_{\chi }$, where the maximum is taking over the
characters of the Hecke torus $T_{A}.$
\end{theorem}

Theorem \ref{QUESS} follows from the fact that the Hecke torus $T_{A}$ acts
on the spaces $\mathcal{H}_{\lambda },$ hence, we can use the Hecke
eigenstates and the bound (\ref{simplified3}). \ 

In particular, using the bound (\ref{BQUESS}), and the explicit information
on $m$ (see Theorem \ref{MF}) we obtain:

\begin{theorem}[Quantum unique ergodicity for statistical states]
Consider an observable $f$ $\in C^{\infty }(\mathbb{T)}$ and a sufficiently
large prime number $p.$ Then 
\begin{equation*}
\left \vert Tr(\pi (f)D_{\lambda })-\tint \limits_{\mathbb{T}}fd\mu \right
\vert \leq \frac{C_{f}}{\sqrt{p}^{N}},
\end{equation*}%
where $\mu =|\omega |^{N}$ is the corresponding volume form and $C_{f}$ is
an explicit computable constant which depends only on the function $f.$
\end{theorem}

\subsection{Results}

\begin{enumerate}
\item \textit{Bounds on higher-dimensional exponential sums.} The main
result of this paper is a new method to bound certain higher-dimensional
exponential sums associated with tori in $Sp(2N,\mathbb{F}_{q}).$ Our method
is based on the self-reducibility property of the Weil representation. As an
application we prove the Hecke quantum unique ergodicity theorem for generic
linear symplectomorphisms of the higher-dimensional tori.

\item \textit{Self-reducibility of the Weil representation. }The main
technical result of this paper is the proof of the self-reducibility of the
Weil representation. This property was described first by G\`{e}rardin in 
\cite{Ge}. However, our proof is slightly different and in particular
applies to the Weil representation over any local field of characteristic
different from two. In order to keep the paper self contained we decided to
present our proof in detail.

\item \textit{Multiplicities. }We present a simple method, using the
self-reducibility property, to compute the dimension of the character spaces
for the action of the tori in the Weil representations.

\item \textit{Two-dimensional Wigner distributions. }We describe a new proof
for the bound on the Wigner distributions associated with tori in $SL(2,%
\mathbb{F}_{q})$. It uses direct geometric calculations, using the new
character formula (\ref{Char H-W}), and avoids the use of the equivariant
property of the Deligne sheaf \cite{GH3}.
\end{enumerate}

\subsection{Structure of the paper}

Apart from the introduction, the paper consists of five sections and two
appendices. In Section \ref{H--WR} we present preliminaries from the theory
of the Heisenberg--Weil representation. Section \ref{selfreducibility}
constitutes the main technical part of this work. Here we formulate and
prove the self-reducibility property of the Weil representation. Section \ref%
{B} deals with the main application of the paper---bounds on the
higher-dimensional Wigner distributions. In Section \ref{HBM} we introduce
the Hannay--Berry model of quantum mechanics on the higher-dimensional tori
and in Section \ref{HQUE} we apply the bounds on the Wigner distribution to
obtain the Hecke quantum unique ergodicity theorem. Finally, in Appendices %
\ref{Pr} and \ref{Pbound} we supply the proofs for the statements that
appear in the body of the paper.

\subsection{Acknowledgments}

It is a pleasure to thank our teacher J. Bernstein. In addition, we thank D.
Kazhdan and M. Baruch for interesting discussions. We are grateful to D.
Kelmer for sharing with us computer simulation data. Finally, we would like
to thank O. Ceyhan and the organizers of the conference AGAQ, Istanbul, June
2006, and J. Wolf and the organizers of the conference Lie Groups, Lie
Algebras and Their Representations, Berkeley, November 2006, for the
invitation to present this work.

\section{The Heisenberg--Weil representation\label{H--WR}}

In this section, we denote by $k=\mathbb{F}_{q}$ the finite field of $q$
elements and odd characteristic.

\subsection{The Heisenberg representation}

Let $(V,\omega )$ be a $2N$-dimensional symplectic vector space over the
finite field $k$. There exists a two-step nilpotent group $H=H\left(
V,\omega \right) $ associated to the symplectic vector space $(V,\omega )$.
\ The group $H$ is called the \textit{Heisenberg group.} It can be realized
as the set $H=V\times k$ equipped with the following multiplication rule: 
\begin{equation*}
(v,z)\cdot (v^{\prime },z^{\prime })=(v+v^{\prime },z+z^{\prime }+\tfrac{1}{2%
}\omega (v,v^{\prime })).
\end{equation*}%
The center of $H$ is $Z(H)=\{(0,z):z\in k\}$. $\ $Fix a non-trivial central
character $\psi :Z(H)\longrightarrow 
%TCIMACRO{\U{2102} }%
%BeginExpansion
\mathbb{C}
%EndExpansion
^{\ast }$. We have the following fundamental theorem:

\begin{theorem}[Stone--von Neumann]
\label{SVN}

There exists a unique $($up to isomorphism$)$ irreducible representation $%
(\pi ,H,\mathcal{H})$ with central character $\psi $, i.e., $\pi (z)=\psi
(z)\cdot Id_{\mathcal{H}}$ for every $z\in Z(H).$
\end{theorem}

We call the representation $\pi $ appearing in Theorem \ref{SVN} the \textit{%
Heisenberg representation} associated with the central character $\psi $.

\subsection{The Weil representation\label{WR}}

Let $Sp=Sp(V,\omega )$ denote the group of linear symplectic automorphisms
of $V$. The group $Sp$ acts by group automorphisms on the Heisenberg group
through its action on the vector space $V$, i.e., $g\cdot (v,z)=(gv,z)$. A
direct consequence of Theorem \ref{SVN} is the existence of a projective
representation $\widetilde{\rho }:Sp\rightarrow PGL(\mathcal{H}).$ The
classical construction \cite{W} works as follows. Considering the Heisenberg
representation $\pi $ and an element $g\in Sp$ we define a new
representation $\pi ^{g}$ acting on the same Hilbert space via $\pi
^{g}(h)=\pi (g(h)).$ Because these irreducible representations share the
same central character then by Theorem \ref{SVN} they are isomorphic and $%
\dim $ $Hom_{H}(\pi ,\pi ^{g})=1$. Choosing for every $g\in Sp$ a non-zero
operator $\widetilde{\rho }(g)\in Hom_{H}(\pi ,\pi ^{g})$ we obtained the
required projective representation. In other words the projective
representation $\widetilde{\rho }$ is characterized by the formula 
\begin{equation}
\widetilde{\rho }(g)\pi (h)\widetilde{\rho }(g)^{-1}=\pi (g(h)),
\label{Egorov1}
\end{equation}%
for every $g\in Sp$, $h\in H.$

It is a deep fact that over finite fields of odd characteristic this
projective representation has a linearization that we will call the Weil
representation.

\begin{theorem}[Weil representation]
\label{linearization}There there exists a unique\footnote{%
Unique except in the case when $k=\mathbb{F}_{3}$ and $\dim (V)=2.$ For the
natural choice in this case see \cite{GH7}.} representation 
\begin{equation*}
\rho :Sp\longrightarrow GL(\mathcal{H)},
\end{equation*}%
satisfying the identity $($\ref{Egorov1}$).$
\end{theorem}

\subsection{The Heisenberg--Weil representation}

Let $J$ denote the semi-direct product $J=Sp\ltimes H.$ The group $J$ is
sometimes referred to as the \textit{Jacobi} group. The compatible pair $%
(\rho ,\pi )$ is equivalent to a single representation 
\begin{equation*}
\tau :J\longrightarrow GL(\mathcal{H}),
\end{equation*}%
of the Jacobi group defined by the formula $\tau (g,h)=\rho (g)\pi (h)$. \
In this paper, we would like to adopt the name \textit{Heisenberg--Weil}
representation for referring to the representation $\tau $.

\subsection{The character of the Heisenberg--Weil representation}

The absolute value of the characters $ch_{\rho }:$ $Sp\rightarrow 
%TCIMACRO{\U{2102} }%
%BeginExpansion
\mathbb{C}
%EndExpansion
$ of the Weil representation and $ch_{\tau }:$ $J\rightarrow 
%TCIMACRO{\U{2102} }%
%BeginExpansion
\mathbb{C}
%EndExpansion
$ of the Heisenberg--Weil representation was described in \cite{H1, H2}, but
the phases have been made explicit only recently in \cite{GH6}. Denote by $%
\sigma :\mathbb{F}_{q}^{\ast }\rightarrow 
%TCIMACRO{\U{2102} }%
%BeginExpansion
\mathbb{C}
%EndExpansion
^{\ast }$ the Legendre (quadratic) character. The following formulas are
taken from \cite{GH6}:

\begin{equation}
ch_{\rho }(g)=\sigma ((-1)^{N}\cdot \det (g-I)),  \label{Char W}
\end{equation}%
\begin{equation}
ch_{\tau }(g,v,z)=ch_{\rho }(g)\cdot \psi (\tfrac{1}{2}\omega (\frac{v}{g-I}%
,v)+z),  \label{Char H-W}
\end{equation}%
for every $g\in Sp$ such that $g-I$ is invertible and every ($v,z)\in H.$

\subsection{Application to multiplicities}

Let us start with the two-dimensional case. Let $T\subset Sp\simeq SL(2,%
\mathbb{F}_{q})$ be a maximal torus. The torus $T$ acts semisimply on $%
\mathcal{H}$, decomposing it into a direct sum of character spaces $\mathcal{%
H=}\oplus \mathcal{H}_{\chi }$ over the characters of $T.$ As a consequence
of having the explicit formula (\ref{Char W}), we obtain\ a simple
description for the multiplicities $m_{\chi }=\dim $ $\mathcal{H}_{\chi }$
(cf. \cite{AM, Ge, Si})$.$ Denote by $\sigma _{T}:T\rightarrow 
%TCIMACRO{\U{2102} }%
%BeginExpansion
\mathbb{C}
%EndExpansion
^{\ast }$ the unique quadratic character of $T.$

\begin{theorem}[Multiplicities formula]
\label{Multiplicities}We have $m_{\chi }=1$ for any character $\chi \neq $ $%
\sigma _{T}.$ Moreover, $m_{\sigma _{T}}=2$ or $0,$ depending on whether the
torus $T$ is split or inert, respectively.
\end{theorem}

For a proof see Appendix \ref{PM}.

Using the orthogonality relation for characters we obtain:

\begin{corollary}
\label{equ}The character $ch_{\rho }$ when restricted to the punctured torus 
$T\smallsetminus I\subset Sp$ equals $\sigma _{T}$ or $-\sigma _{T}$
depending on whether $T$ is split or inert, respectively.
\end{corollary}

In Subsection \ref{HDM} we use the self-reducibility property and extend
Theorem \ref{Multiplicities} to the higher-dimensional Weil representations.

\section{Self-reducibility of the Weil representation\label{selfreducibility}%
}

In this section, unless stated otherwise, the field $k$ is an arbitrary
local or finite field of characteristic different from two.

\subsection{Symplectic module structures\label{Sms}}

Let $K$ be a finite-dimensional commutative algebra over the field $k$. \
Let $Tr:$ $K\rightarrow k$ be the trace map associating to an element $x\in $
$K$ the trace of the $k$-linear operator $m_{x}:$ $K\rightarrow $ $K$
obtained by left multiplication by the element $x$. Consider a symplectic
vector space $(V,\omega )$ over $k.$

\begin{definition}
\label{SMS}A symplectic $K$-module structure on $\left( V,\omega \right) $
is an action $K\otimes _{k}V\rightarrow V,$ and a $K$-linear symplectic form 
$\overline{\omega }:V\times V\rightarrow $ $K$ such that 
\begin{equation}
Tr\circ \overline{\omega }=\omega .  \label{SC}
\end{equation}
\end{definition}

Given a symplectic module structure $(K,V,\overline{\omega })$ on a
symplectic vector space $(V,\omega ),$ we denote by $\overline{Sp}=Sp(V,%
\overline{\omega })$ the group of $K$-linear symplectomorphisms with respect
to the form $\overline{\omega }$. The compatibility condition (\ref{SC})
gives a natural \ embedding 
\begin{equation}
\iota _{S}:\overline{Sp}\hookrightarrow Sp.  \label{im}
\end{equation}%
.

\subsection{Symplectic module structure associated with a maximal torus\label%
{SmsT}}

Let $T\subset Sp$ be a maximal torus.

\subsubsection{A \textbf{particular case\label{particular}}}

For simplicity, let us assume first that $T$ acts irreducibly on the vector
space $V$, i.e., there exists no non-trivial $T$-invariant subspaces. Let $%
A=Z(T,End(V))$\ be the centralizer of $T$ in the algebra of all linear
endomorphisms. Clearly (due to the assumption of irreducibility) $A$ is a
division algebra. In addition we have:\  \  \  \  \  \  \  \  \  \  \  \  \  \  \  \  \  \  \
\  \  \  \  \  \  \  \  \  \  \  \  \  \  \  \  \  \  \  \  \  \  \  \  \  \  \  \  \  \  \  \  \  \  \  \  \  \
\  \  \  \  \  \  \  \  \  \  \  \  \  \  \  \  \  \  \  \  \  \  \  \  \  \  \  \  \  \  \  \  \  \  \  \  \  \
\  \  \  \  \  \  \  \  \  \  \  \  \  \  \  \  \  \  \  \  \  \  \  \  \  \  \  \  \  \  \  \  \  \  \  \  \  \
\  \  \  \  \  \  \  \  \  \  \  \  \  \  \  \  \  \  \  \  \  \  \  \  \  \  \  \  \  \  \  \  \  \  \  \  \  \
\  \  \  \  \  \  \  \  \  \  \  \  \  \  \  \  \  \  \  \  \  \  \  \  \  \  \  \  \  \  \  \  \  \  \  \  \  \
\  \  \  \  \  \  \  \  \  \  \  \  \  \  \  \  \  \  \  \  \  \  \  \  \  \  \  \  \  \  \  \  \  \  \  \  \  \
\  \  \  \  \  \  \  \  \  \  \  \  \  \  \  \  \  \  \  \  \  \  \  \  \  \  \  \  \  \  \  \  \  \  \  \  \  \
\  \  \  \  \  \  \  \  \  \  \  \  \  \  \  \  \  \  \  \  \  \  \  \  \  \  \  \  \  \  \  \  \  \  \  \  \  \
\  \  \  \  \  \  \  \  \  \  \  \  \  \  \  \  \  \  \  \  \  \  \  \  \  \  \  \  \  \  \  \  \  \  \  \  \  \
\  \  \  \  \  \  \  \  \  \  \  \  \  \  \  \  \  \  \  \  \  \  \  \  \  \  \  \  \  \  \  \  \  \  \  \  \  \
\  \  \  \  \  \  \  \  \  \  \  \  \  \  \  \  \  \  \  \  \  \  \  \  \  \  \  \  \  \  \  \  \  \  \  \  \  \
\  \  \  \  \  \  \  \  \  \  \  \  \  \  \  \  \  \  \  \  \  \  \  \  \  \  \  \  \  \  \  \  \  \  \  \  \  \
\  \  \  \  \  \  \  \  \  \  \  \  \  \  \  \  \  \  \  \  \ 

\begin{claim}
\label{Comm}The algebra $A$ is commutative.
\end{claim}

For a proof see Appendix \ref{ComDim}.

In particular, Claim \ref{Comm} implies that $A$ is a field. Let us now
describe a special quadratic element in the Galois group $Gal(A/k)$ of all
the automorphisms of $A$ which leave the field $k$ fixed. Denote by $(\cdot
)^{t\text{ }}:End(V)\rightarrow End(V)$ the symplectic transpose
characterized by the property $\omega (Rv,u)=\omega (v,R^{t}u)$ for all $%
v,u\in V$, and every $R\in End(V)$. It can be easily verified that $(\cdot
)^{t}$ preserves $A$, leaving the subfield $k$ fixed, hence, it defines an
element $\Theta \in Gal(A/k)$ satisfying $\Theta ^{2}=Id$\text{. \ Denote by 
}%
\begin{equation*}
K=A^{\Theta },
\end{equation*}%
\text{the subfield of }$A$\text{ consisting of the elements fixed by }$%
\Theta $\text{. }

\begin{proposition}[Hilbert's theorem 90]
\label{prop_dim}We have $\dim _{K}V=2.$
\end{proposition}

\text{For a proof see Appendix \ref{ComDim}.}

\begin{corollary}
We have\ $\dim _{K}A=2.$
\end{corollary}

As a corollary, we have the following description of $T$. Denote by $%
N_{A/K}:A\rightarrow K$ the standard norm map.

\begin{corollary}
\label{DT}We have $T=S(A)=\left \{ a\in A:N_{A/K}(a)=1\right \} .$
\end{corollary}

For a proof see Appendix \ref{Norm}.

The symplectic form $\omega $ can be lifted to a $K$-linear symplectic form $%
\overline{\omega }$ which is invariant under the action of the torus $T$.
This is the content of the following proposition:

\begin{proposition}[Existence of canonical symplectic module structure]
\label{canonicallift}There exists a\ canonical $T$-invariant $K$-linear
symplectic form $\overline{\omega }:V\times V\rightarrow K$ satisfying the
property $Tr\circ \overline{\omega }=\omega .$
\end{proposition}

For a proof see Appendix \ref{C-L}.

Concluding, we obtained a $T$-invariant symplectic $K$-module structure on $%
V $.

Let $\overline{Sp}=Sp(V,\overline{\omega })$ denote the group of $K$-linear
symplectomorphisms with respect to the symplectic form $\overline{\omega }$.
We denote the embedding (\ref{im}) by $i_{S}:\overline{Sp}\hookrightarrow
Sp. $ The elements of $T$ commute with the action of $K$ and preserve the
symplectic form $\overline{\omega }$\ (Proposition \ref{canonicallift});
hence, we can consider $T$ as a subgroup of $\overline{Sp}$. \ By
Proposition \ref{prop_dim} we can identify $\overline{Sp}\simeq SL(2,K),$
and using (\ref{im}) we obtain 
\begin{equation}
T\subset SL(2,K)\subset Sp.  \label{inr}
\end{equation}%
Moreover, we see that $T$ consists of the set of $K$-rational points of a
maximal algebraic torus $\mathbf{T\subset SL}_{2}$ (we use bold face to
denote algebraic varieties).

\subsubsection{\textbf{General case\label{general}}}

Here, we drop the assumption that $T$ acts irreducibly on $V$. By the same
argument as before one can show that the algebra $A=Z(T,End(V))$ is
commutative, yet, it may no longer be a field. The symplectic transpose $%
(\cdot )^{t}$ preserves the algebra $A$, and induces an involution $\Theta
:A\rightarrow A.$ Let $K=A^{\Theta }$ be the subalgebra consisting of
elements $a\in A$ fixed by $\Theta $. \ Following the same argument as in
the proof of Proposition \ref{prop_dim}, one shows that $V$ is a free $K$%
-module of rank $2$. \ Following the same arguments as in the proof of
Proposition \ref{canonicallift}, one shows that there exists a canonical
symplectic form $\overline{\omega }:V\times V\rightarrow K$, which is $K$%
-linear and invariant under the action of the torus $T$. Concluding,
associated to a maximal torus $T$ there exists a $T$-invariant symplectic $K$%
-module structure\ $(K,V,\overline{\omega }).$

Denote by $\overline{Sp}=Sp(V,\overline{\omega })$ the group of $K$-linear
symplectomorphisms with respect to the form $\overline{\omega }.$ We have a
natural embedding $\iota _{S}:\overline{Sp}\hookrightarrow Sp$ and we can
consider $T$ as a subgroup of $\overline{Sp}$. Finally, we have $\overline{Sp%
}\simeq SL(2,K)$ and $T$ consists of the $K$-rational points of a maximal
torus $\mathbf{T\subset SL}_{2}$. In particular, the relation (\ref{inr})
holds also in this case%
\begin{equation*}
T\subset SL(2,K)\subset Sp.
\end{equation*}

\ We shall now proceed to give a finer description of all objects discussed
so far. The main technical result is summarized in the following lemma:

\begin{lemma}[Symplectic decomposition]
\label{technical}We have a canonical decomposition 
\begin{equation}
(V,\omega )=\tbigoplus \limits_{\alpha \in \Xi }(V_{\alpha },\omega _{\alpha
}),  \label{SR}
\end{equation}%
into $(T,A)$-invariant symplectic subspaces. In addition, we have the
following associated canonical decompositions:
\end{lemma}

\begin{enumerate}
\item $T=\tprod T_{\alpha },$ where $T_{\alpha }$ consists of elements $g\in
T$ such that $g_{|V_{\beta }}=Id$ for every $\beta \neq \alpha $.

\item $A=\oplus A_{\alpha }$, where $A_{\alpha }$ consists of elements $a\in
A$ such that $a_{|V_{\beta }}=Id$ for every $\beta \neq \alpha $. Moreover,
each subalgebra $A_{\alpha }$ is preserved under the involution $\Theta .$

\item $K=\oplus K_{\alpha }$, where $K_{\alpha }=A_{\alpha }^{\Theta }$.
Moreover, $K_{\alpha }$ is a field and $\dim _{K_{\alpha }}V_{\alpha }=2$.

\item $\overline{\omega }=\oplus \overline{\omega }_{\alpha }$, where $%
\overline{\omega }_{\alpha }:V_{\alpha }\times V_{\alpha }\rightarrow
K_{\alpha }$ is a $K_{\alpha }$-linear $T_{\alpha }$-invariant symplectic
form satisfying $Tr\circ \overline{\omega }_{\alpha }=\omega _{\alpha }$.
\end{enumerate}

For a proof see Appendix \ref{Ptechnical}.

\begin{definition}
\label{type}We will call the set $\Xi $\ $\left( \text{\ref{SR}}\right) $
the symplectic type of $T$ and the number $\left \vert \Xi \right \vert $
the symplectic rank of $T$.
\end{definition}

Using the results of Lemma \ref{technical}, we have an isomorphism 
\begin{equation}
\overline{Sp}\simeq \tprod \overline{Sp}_{\alpha },  \label{decomp5}
\end{equation}%
where $\overline{Sp}_{\alpha }=Sp(V_{\alpha },\overline{\omega }_{\alpha })$
denotes the group of $K_{\alpha }$-linear symplectomorphisms with respect to
the form $\overline{\omega }_{\alpha }$. Moreover, \ for every $\alpha \in
\Xi $ we have $T_{\alpha }\subset \overline{Sp}_{\alpha }$. In particular,
under the identifications $\overline{Sp}_{\alpha }\simeq SL(2,K_{\alpha }),$
there exist the following sequence of inclusions of groups:%
\begin{equation}
T=\tprod T_{\alpha }\subset \tprod SL(2,K_{\alpha })=SL(2,K)\subset Sp,
\label{IOG}
\end{equation}%
and for every $\alpha \in \Xi $ the torus $T_{\alpha }$ coincides with the $%
K_{\alpha }$-rational points of a maximal torus $\mathbf{T}_{\alpha }\subset 
\mathbf{SL}_{2}$.

\subsection{Self-reducibility of the Weil representation\label{selfred}}

In this subsection we assume that the field $k$ is a finite field\footnote{%
We remark that the results continue to hold true also for local fields of
characteristic $\neq 2$, i.e., with the appropriate modification, replacing
the group $Sp$ with its double cover $\widetilde{Sp}$ \cite{W}.}. Let $(\tau
,J,\mathcal{H)}$ be the Heisenberg--Weil representation associated with a
central character $\psi :Z(J)=Z(H)\rightarrow 
%TCIMACRO{\U{2102} }%
%BeginExpansion
\mathbb{C}
%EndExpansion
^{\ast }$. \ Recall that $J=Sp\ltimes H$ and $\tau $ is obtained as a
semi-direct product $\tau =\rho \ltimes \pi $ of the Weil representation $%
\rho $ and the Heisenberg representation $\pi $. Let $T\subset Sp$ be a
maximal torus.

\subsubsection{A particular case\textbf{\label{particular1}}}

For clarity of presentation, let us assume first that $T$ acts irreducibly
on $V$. Using the results of the previous section, there exists a symplectic
module structure $(K,V,\overline{\omega })$ where $K/k$ is a field extension
of degree $[K:k]=N$. The group $\overline{Sp}=Sp(V,\overline{\omega })$ is
embedded as a subgroup $\iota _{S}:\overline{Sp}\hookrightarrow Sp.$ Our
goal is to describe the restriction 
\begin{equation}
(\overline{\rho }=\iota _{S}^{\ast }\rho ,\overline{Sp},\mathcal{H)}.
\label{RR}
\end{equation}%
Define an auxiliary Heisenberg group 
\begin{equation}
\overline{H}=V\times K,  \label{auxH}
\end{equation}%
and the multiplication is given by $(v,z)\cdot (v^{\prime },z^{\prime
})=(v+v^{\prime },z+z^{\prime }+%
%TCIMACRO{\U{bd}}%
%BeginExpansion
{\frac12}%
%EndExpansion
\overline{\omega }(v,v^{\prime }))$.

There exists a homomorphism%
\begin{equation}
\iota _{H}:\overline{H}\rightarrow H,  \label{homauxH}
\end{equation}%
given by $(v,z)\mapsto (v,Tr(z))$. \ 

Consider the pullback $(\overline{\pi }=\iota _{H}^{\ast }\pi ,\overline{H},%
\mathcal{H)}.$

\begin{proposition}
\label{resHeisenberg}The representation $(\overline{\pi }=$ $\iota
_{H}^{\ast }\pi ,\overline{H},\mathcal{H)}$ is the Heisenberg representation
associated with the central character $\overline{\psi }=\psi \circ Tr$.
\end{proposition}

For a proof see Appendix \ref{PresHeisenberg}.

The group $\overline{Sp}$ acts by automorphisms on the group $\overline{H}$
through its tautological action on the $V$-coordinate. This action is
compatible with the action of $Sp$ on $H$, i.e., we have $\iota _{H}(g\cdot
h)=$ $\iota _{S}(g)\cdot $ $\iota _{H}(h)$ for every $g\in \overline{Sp}$,
and $h\in \overline{H}$. \ The description of the representation $\overline{%
\rho }$ \ (\ref{RR}) now follows easily (cf. \cite{Ge})

\begin{theorem}[Self-reducibility property---particular case]
\label{reWeil}The representation $(\overline{\rho },\overline{Sp},\mathcal{H)%
}$ is the Weil representation associated with the Heisenberg representation $%
(\overline{\pi },\overline{H},\mathcal{H)}$.
\end{theorem}

For a proof see Appendix \ref{PreWeil}.

\begin{remark}
We can summarize the result in a slightly more elegant manner using the
Jacobi groups. Let $J=Sp\ltimes H$ and $\overline{J}=\overline{Sp}\ltimes 
\overline{H}$ be the Jacobi groups associated with the symplectic spaces $%
(V,\omega )$ and $(V,\overline{\omega })$ respectively. We have a
homomorphism $\iota :\overline{J}\rightarrow J$ given by $\iota (g,h)=($ $%
\iota _{S}(g),$ $\iota _{H}(h))$. Let $(\tau ,J,\mathcal{H)}$ be the
Heisenberg--Weil representation of $J$ associated with a character $\psi $
of the center $Z(J)$ (note that $Z(J)=Z(H)$), then the pullback $(\iota
^{\ast }\tau ,\overline{J},\mathcal{H)}$ is the Heisenberg--Weil
representation of $\overline{J}$ associated with the character $\overline{%
\psi }=\psi \circ Tr$ of the center $Z(\overline{J})$.
\end{remark}

\subsubsection{The general case\textbf{\label{general1}}}

Here, we drop the assumption that $T$ acts irreducibly on $V$. Let $(K,V,%
\overline{\omega })$ be the associated symplectic module structure. Using
the results of Subsection \ref{general}, we have decompositions 
\begin{equation}
(V,\omega )=\tbigoplus \limits_{\alpha \in \Xi }(V_{\alpha },\omega _{\alpha
}),\text{ \ }(V,\overline{\omega })=\tbigoplus \limits_{\alpha \in \Xi
}(V_{\alpha },\overline{\omega }_{\alpha }),  \label{decomm7}
\end{equation}%
where $\overline{\omega }_{\alpha }:V_{\alpha }\times V_{\alpha }\rightarrow
K_{\alpha }$. Let $\overline{H}=V\times K$ be the Heisenberg group
associated with $(V,\overline{\omega })$ (cf. \ref{auxH}). There exists (cf. %
\ref{homauxH}) an homomorphism $\iota _{H}:\overline{H}\rightarrow H.$ Let
us describe the pullback $\overline{\pi }=\iota _{H}^{\ast }\pi $ of the
Heisenberg representation. First, we note that the decomposition (\ref%
{decomm7}) induces a corresponding decomposition of the Heisenberg group $%
\overline{H}=\tprod \overline{H}_{\alpha },$ where $\overline{H}_{\alpha }$
is the Heisenberg group associated with $(V_{\alpha },\overline{\omega }%
_{\alpha })$.

\begin{proposition}
\label{Heisenbergdecom}There exists an isomorphism 
\begin{equation*}
(\overline{\pi },\overline{H},\mathcal{H)\simeq (}\tbigotimes \overline{\pi }%
_{\alpha },\tprod \overline{H}_{\alpha },\tbigotimes \mathcal{H}_{\alpha }),
\end{equation*}%
where $(\overline{\pi }_{\alpha },\overline{H}_{\alpha },\mathcal{H}_{\alpha
})$ is the Heisenberg representation of $\overline{H}_{\alpha }$ associated
with the central character $\overline{\psi }_{\alpha }=\psi \circ
Tr_{K_{\alpha }/k}$.
\end{proposition}

For a proof see Appendix \ref{PHeisenbergdecom}.

Let $\iota _{S}:\overline{Sp}\hookrightarrow Sp$ be the embedding (\ref{im}%
). Our next goal is to describe (cf. \cite{Ge}) the restriction $\overline{%
\rho }=\iota _{S}^{\ast }\rho .$ Recall the decomposition $\overline{Sp}%
=\tprod \overline{Sp}_{\alpha }$ (see (\ref{decomp5})).

\begin{theorem}[Self-reducibility property---general case]
\label{Weildecomp}There exists an isomorphism 
\begin{equation*}
(\overline{\rho },\overline{Sp},\mathcal{H)\simeq (}\tbigotimes \overline{%
\rho }_{\alpha },\tprod \overline{Sp}_{\alpha },\tbigotimes \mathcal{H}%
_{\alpha }),
\end{equation*}%
where $(\overline{\rho }_{\alpha },\overline{Sp}_{\alpha },\mathcal{H}%
_{\alpha })$ is the Weil representation associated with the Heisenberg
representation\ $\overline{\pi }_{\alpha }$.
\end{theorem}

For a proof see Appendix \ref{PWeildecomp}.

\begin{remark}
As before, we can state an equivalent result using the Jacobi groups $%
J=Sp\ltimes H$ and $\overline{J}=\overline{Sp}\ltimes \overline{H}$. We have
a decomposition $\overline{J}=\tprod \overline{J}_{\alpha },$ where $%
\overline{J}_{\alpha }=\overline{Sp}_{\alpha }\ltimes \overline{H}_{\alpha }$%
. Let $\tau $ be the Heisenberg--Weil representation of $J$ associated with
a character $\psi $ of the center $Z(J)$ $($note that $Z(J)=Z(H))$. Then the
pullback $\overline{\tau }=$ $\iota ^{\ast }\tau $ is isomorphic to $\otimes 
\overline{\tau }_{\alpha }$, where $\overline{\tau }_{\alpha }$ is the
Heisenberg--Weil representation of $\overline{J}_{\alpha }$, associated with
the character $\overline{\psi }_{\alpha }=\psi \circ Tr_{K_{\alpha }/k}$ of
the center $Z(\overline{J}_{\alpha })$.
\end{remark}

\subsection{Application to multiplicities\label{HDM}}

Let $T$ $\subset Sp$ be a maximal torus. The torus $T$ acts, via the Weil
representation $\rho $, on the space $\mathcal{H}$, decomposing it into a
direct sum of character spaces $\mathcal{H=}\oplus \mathcal{H}_{\chi }$. We
would like to compute the multiplicities $m_{\chi }=\dim (\mathcal{H}_{\chi
})$. Using Lemma \ref{technical}, we have (see (\ref{IOG})) a canonical
decomposition of $T\ $%
\begin{equation}
T=\tprod T_{\alpha },  \label{dec}
\end{equation}%
where each of the tori $T_{\alpha }$ coincides with a maximal torus inside $%
\overline{Sp}\simeq SL(2,K_{\alpha }),$ for some field extension $K_{\alpha
}\supset k.$ In particular, by (\ref{dec}) we have a decomposition 
\begin{equation}
\mathcal{H}_{\chi }\mathcal{=}\tbigotimes \limits_{\chi _{\alpha }:T_{\alpha
}\rightarrow 
%TCIMACRO{\U{2102} }%
%BeginExpansion
\mathbb{C}
%EndExpansion
^{\ast }}\mathcal{H}_{\chi _{\alpha }},  \label{decchar}
\end{equation}%
where $\chi =\tprod \chi _{\alpha }:\tprod T_{\alpha }\rightarrow 
%TCIMACRO{\U{2102} }%
%BeginExpansion
\mathbb{C}
%EndExpansion
^{\ast }$. Hence, by Theorem \ref{Weildecomp} and the result about the
multiplicities in the two-dimensional case (see Proposition (\ref%
{Multiplicities})), we can compute the integer $m_{\chi }.$ Denote by $%
\sigma _{T_{\alpha }\text{ }}$the quadratic character of $T_{\alpha }$ (note
that by Proposition (\ref{Multiplicities}) the quadratic character $\sigma
_{T_{\alpha }}$ cannot appear in the decomposition (\ref{decchar}) if the
torus $T_{\alpha }$ is inert).

\begin{theorem}[Multiplicities formula---higher dimensional]
\label{MF}We have 
\begin{equation*}
m_{\chi }=2^{l},
\end{equation*}%
where $l=|\{ \alpha :$ $\chi _{\alpha }=$ $\sigma _{T_{\alpha }}\}|.$
\end{theorem}

\section{Bounds on higher-dimensional Wigner distributions\label{B}}

Let $(\rho ,Sp,\mathcal{H)}$ be the Weil representation associated with a $%
2N $-dimensional vector space over an odd characteristic finite field $k=%
\mathbb{F}_{q}$. Consider a maximal torus $T\subset Sp$ and the associated
decomposition of $\mathcal{H}$ into a direct sum of character spaces $%
\mathcal{H=}\oplus \mathcal{H}_{\chi }.$ For a character vector $\varphi \in 
\mathcal{H}_{\chi }$ we will bound the Wigner distribution $\left \langle
\varphi |\pi (v)\varphi \right \rangle $ where $v\in V$ is not contained in
any proper $T$-invariant subspace. We will explain how to use the
self-reducibility property for this purpose.

\subsection{The completely inert case}

It will be convenient to assume first that the torus $T$ is completely
inert, i.e., acts irreducibly on $V$.

\begin{theorem}[Bound on Wigner distributions---inert case]
\label{KR1}For every non-zero vector $v\in V$ we have 
\begin{equation*}
\left \vert \left \langle \varphi |\pi (v)\varphi \right \rangle \right
\vert \leq \frac{2}{\sqrt{q}^{N}}.
\end{equation*}
\end{theorem}

To get this bound we proceed as follow. The torus $T$ acts irreducibly on
the vector space $V$. Invoking the result of Section \ref{particular}, there
exists a canonical symplectic module structure $(K,V,\overline{\omega })$
associated to $T$. Recall that in this particular case the algebra $K$ is in
fact a field, $K=\mathbb{F}_{q^{N}},$ and $\dim _{K}V=2$. Let $\overline{J}=%
\overline{Sp}\ltimes \overline{H}$ be the Jacobi group associated to the
two-dimensional symplectic vector space $(V,\overline{\omega })$. There
exists a natural homomorphism $\iota :\overline{J}\rightarrow J$. \ Invoking
the results of\ Section \ref{particular1}, the pullback $\overline{\tau }=$ $%
\iota ^{\ast }\tau $ is the Heisenberg--Weil representation of $\overline{J}$%
, i.e., $\overline{\tau }=\overline{\rho }\ltimes \overline{\pi }$. \ The
orthogonal projector $P_{\chi }$ on the space $\mathcal{H}_{\chi }$ can be
written in terms of the Weil representation $\overline{\rho }$ as $P_{\chi
}=|T|^{-1}\sum_{g\in T}\chi ^{-1}(g)\overline{\rho }(g).$ Since $\dim 
\mathcal{H}_{\chi }=1$ (Theorem \ref{MF}) we realize that $\left \langle
\varphi |\pi (v)\varphi \right \rangle =Tr(P_{\chi }\pi (v))$ where $\varphi
\in \mathcal{H}_{\chi }.$ Overall, we have 
\begin{equation}
\left \langle \varphi |\pi (v)\varphi \right \rangle =\frac{1}{|T|}\tsum
\limits_{g\in T}\chi ^{-1}(g)Tr(\overline{\rho }(g)\pi (v)).  \label{Wig}
\end{equation}%
Note that $Tr(\overline{\rho }(g)\pi (v))$ is nothing other than the
character $ch_{\overline{\tau }}(g\cdot v)$ of the Heisenberg--Weil
representation $\overline{\tau }$ and that $|T|=q^{N}+1.$ Therefore the
right-hand side of (\ref{Wig}) is defined completely in terms of the
two-dimension\textbf{\ }Heisenberg--Weil representation $\overline{\tau }$.
\ Theorem \ref{KR1} is then a particular case of the following theorem:

\begin{theorem}[\protect \cite{GH3}]
\label{twodim1}Let $(V,\omega )$ be a two-dimensional symplectic vector
space over the finite field $k=\mathbb{F}_{q}$, and $(\tau ,J,\mathcal{H)}$
the corresponding Heisenberg--Weil representation. Let $T\subset Sp$ be a
maximal torus. We have the following estimate: 
\begin{equation}
\left \vert \tsum \limits_{g\in T}\chi (g)ch_{\tau }(g\cdot v)\right \vert
\leq 2\sqrt{q},  \label{2dim_estimate1}
\end{equation}%
where $\chi $ is a character of $T$, and $0\neq v\in V$ is not an
eigenvector of $T$.
\end{theorem}

For the sake of completeness we write in Appendix \ref{Pbound} a new proof
of Theorem \ref{twodim1}.

\begin{remark}[Alternative approach]
We propose another approach to obtain the relation between the higher
dimensional exponential sum over $k=\mathbb{F}_{q}$ and the one dimensional
sum over $K=\mathbb{F}_{q^{N}}.$ We use the existence of symplectic module
structure (Theorem \ref{canonicallift}) and then---instead of using the
self-reducibility property---we invoke the explicit formula (\ref{Char H-W})
for the character of the Heisenberg--Weil representation. Denote by $(\tau
,J,\mathcal{H)}$ the Heisenberg--Weil representation associated with a $2N$%
-dimensional symplectic vectors space over $k$ and by $(K,V,\overline{\omega 
})$ the symplectic module structure associated with the torus $T.$ We deduce
that%
\begin{eqnarray*}
\sum_{g\in T}\chi (g)ch_{\tau }(g\cdot v) &=&\sum_{g\in T}\chi (g)\sigma
((-1)^{N}\det (g-I))\psi (\tfrac{1}{2}\omega (\tfrac{v}{g-I}{\small ,v})) \\
&=&\sum_{g\in T}\chi (g)\overline{\sigma }(-\det {}_{K}(g-I))\overline{\psi }%
(\tfrac{1}{2}\overline{\omega }(\tfrac{v}{g-I}{\small ,v})),
\end{eqnarray*}%
where $\overline{\sigma }=\sigma \circ Norm_{K/k}$, $\overline{\psi }=\psi
\circ Tr_{K/k},$ and $\det_{K}(g-I)$ is the determinant of $g-I$ acting on $%
V $ as a vector space over $K.$ The last sum above is over a torus in \ $%
Sp(V,\overline{\omega })$ which we evaluate using Theorem \ref{twodim1}.
\end{remark}

\subsection{General case}

In this subsection we state and prove the analogue of Theorem \ref{KR1},
where we drop the assumption of $T$ being completely inert. In what follows,
we use the results of Subsections \ref{general} and \ref{general1}.

Let $(K,V,\overline{\omega })$ be the symplectic module structure associated
with the torus $T$. The algebra $K$ is no longer a field, but decomposes
into a direct sum of fields $K=\oplus _{\alpha \in \Xi }K_{\alpha }$. We
have canonical decompositions $(V,\omega )=\oplus (V_{\alpha },\omega
_{\alpha })$ and $(V,\overline{\omega })=\oplus (V_{\alpha },\overline{%
\omega }_{\alpha }).$ Recall that $V_{\alpha }$ is a two-dimensional vector
space over the field $K_{\alpha }$. The Jacobi group $\overline{J}$
decomposes into $\overline{J}=\Pi \overline{J}_{\alpha },$ where $\overline{J%
}_{\alpha }=\overline{Sp}_{\alpha }\ltimes \overline{H}_{\alpha }$ is the
Jacobi group associated to $(V_{\alpha },\overline{\omega }_{\alpha }).$ The
pullback $(\overline{\tau }=$ $\iota ^{\ast }\tau ,\overline{J},\mathcal{H)}$
decomposes into a tensor product $(\otimes \overline{\tau }_{\alpha },\Pi 
\overline{J}_{\alpha },\otimes \mathcal{H}_{\alpha })$ where $\overline{\tau 
}_{\alpha }$ is the Heisenberg--Weil representation of $\overline{J}_{\alpha
}$. The torus $T$ decomposes into $T=\Pi T\alpha $ where $T_{\alpha }$ is a
maximal torus in $\overline{Sp}_{\alpha }.$ Consequently, the character $%
\chi :T\rightarrow 
%TCIMACRO{\U{2102} }%
%BeginExpansion
\mathbb{C}
%EndExpansion
^{\ast }$ decomposes into a product $\chi =\Pi \chi _{\alpha }:\Pi T_{\alpha
}\rightarrow 
%TCIMACRO{\U{2102} }%
%BeginExpansion
\mathbb{C}
%EndExpansion
^{\ast }$ and the space $\mathcal{H}_{\chi }$ decomposes into a tensor
product over the character $\chi _{\alpha }:T_{\alpha }\rightarrow 
%TCIMACRO{\U{2102} }%
%BeginExpansion
\mathbb{C}
%EndExpansion
^{\ast }.$

\begin{equation}
\mathcal{H}_{\chi }=\tbigotimes \mathcal{H}_{\chi _{\alpha }}.  \label{ESD}
\end{equation}

It follows from the above decomposition that it is enough to estimate matrix
coefficients with respect to \textquotedblleft pure tensor\textquotedblright
\ character vector $\varphi $ of the form $\varphi =\otimes \varphi _{\alpha
}$, where $\varphi _{\alpha }\in \mathcal{H}_{\chi _{\alpha }}$. For a
vector of the form $v=\otimes v_{\alpha }$ we have 
\begin{equation}
\left \langle \tbigotimes \varphi _{\alpha }|\pi (v)\tbigotimes \varphi
_{\alpha }\right \rangle =\tprod \left \langle \varphi _{\alpha }|\pi
(v_{\alpha })\varphi _{\alpha }\right \rangle .  \label{MultipF}
\end{equation}

Hence, we are reduced to estimate the matrix coefficients $\left \langle
\varphi _{\alpha }|\pi (v_{\alpha })\varphi _{\alpha }\right \rangle $, but
these are defined in terms of the two-dimensional Heisenberg--Weil
representation $\overline{\tau }_{\alpha }$. In addition, we recall the
assumption that the vector $v\in V$ is not contained in any proper $T$%
-invariant subspace. This condition in turn implies that no summand $%
v_{\alpha }$ \ is an eigenvector of \ $T_{\alpha }.$ Hence, we can use
Theorem \ref{twodim1} and the fact that $|T_{\alpha }|$ is of order of $%
q^{[K_{\alpha }:\mathbb{F}_{q}]}\pm 1$ to get 
\begin{equation}
\left \vert \left \langle \varphi _{\alpha }|\pi (v_{\alpha })\varphi
_{\alpha }\right \rangle \right \vert \leq (2+o(1))/\sqrt{q}^{[K_{\alpha }:%
\mathbb{F}_{q}]}.  \label{Bound}
\end{equation}%
\ Consequently, using $\left( \text{\ref{MultipF}}\right) $ and (\ref{Bound}%
) we obtain 
\begin{equation*}
\left \vert \left \langle \tbigotimes \varphi _{\alpha }|\pi (v)\tbigotimes
\varphi _{\alpha }\right \rangle \right \vert \leq (2+o(1))^{\left \vert \Xi
\right \vert }/\sqrt{q}^{\tsum [K_{\alpha }:\mathbb{F}_{q}]}=(2+o(1))^{\left%
\vert \Xi \right \vert }/\sqrt{q}^{N},
\end{equation*}%
where $r_{p}=\left \vert \Xi \right \vert $ is the symplectic rank of the
torus $T$. Let us summarize:

\begin{theorem}[Bound on Wigner distributions---general case]
\label{KR_general1}Let $(V,\omega )$ be a $2N$-dimensional vector space over
the finite field $\mathbb{F}_{q}$ and $(\tau ,J,\mathcal{H)}$ the
corresponding Heisenberg--Weil representation. Let $\varphi \in \mathcal{H}%
_{\chi }$ be a unit $\chi $-eigenstate with respect to a maximal torus $%
T\subset Sp$. We have the following estimate: 
\begin{equation*}
\left \vert \left \langle \varphi |\pi (v)\varphi \right \rangle \right
\vert \leq \frac{m_{\chi }\cdot (2+o(1))^{r_{p}}}{\sqrt{q}^{N}},
\end{equation*}%
where $1\leq r_{p}\leq N$ is the symplectic rank of $T$, $m_{\chi }=\dim 
\mathcal{H}_{\chi }$, and $v\in V$ is not contained in any $T$-invariant
subspace.
\end{theorem}

\section{The Hannay--Berry model\label{HBM}}

We shall proceed to describe the higher-dimensional Hannay--Berry model of
quantum mechanics on toral phase spaces. This model plays an important role
in the mathematical theory of quantum chaos as it serves as a model where
general phenomena, which are otherwise treated only on a heuristic basis,
can be rigorously proven.

\subsection{The phase space}

Our phase space is the $2N$-dimensional symplectic torus $(\mathbb{T},\omega
).$ We denote by $\Gamma $ the group of linear symplectomorphisms of $%
\mathbb{T}$. Note that $\Gamma \simeq Sp(2N,%
%TCIMACRO{\U{2124} }%
%BeginExpansion
\mathbb{Z}
%EndExpansion
).$ On the torus $\mathbb{T}$ we consider the algebra $\mathcal{A}$ of
observables---trigonometric polynomials. The algebra $\mathcal{A}$ has as a
natural basis the lattice $\Lambda $ of characters (exponents) of $\mathbb{T}
$. The form $\omega $ induces a skew-symmetric form on $\Lambda $, which we
denote also by $\omega $, and we assume it takes integral values\ on $%
\Lambda $ and is normalized so that $\int_{\mathbb{T}}|\omega |^{N}=1.$

\subsection{The mechanical system}

Our mechanical system is of a very simple nature, consists of an
automorphism $A\in \Gamma $ which we assume to be generic element (see
Definition \ref{qergodic}), i.e., $A$ is regular and admits no invariant
co-isotropic sub-tori. The last condition can be equivalently restated in
dual terms, namely, requiring that $A$ admits no invariant isotropic
subvectorspaces in $\Lambda _{%
%TCIMACRO{\U{211a} }%
%BeginExpansion
\mathbb{Q}
%EndExpansion
}=\Lambda \otimes _{%
%TCIMACRO{\U{2124} }%
%BeginExpansion
\mathbb{Z}
%EndExpansion
}%
%TCIMACRO{\U{211a} }%
%BeginExpansion
\mathbb{Q}
%EndExpansion
$. The element $A$ generates, via its action as an automorphism $A:\mathbb{%
T\longrightarrow T},$\ a discrete time dynamical system.

\subsection{Quantization via the non-commutative torus model\label{RM}}

In this paper we employ a quantization model, that we call the \textit{%
non-commutative torus }model, developed in \cite{GH2, GH3, GH9}. This is a
certain one-parameter family of \textquotedblleft
protocols\textquotedblright \ parameterized by a parameter $\hbar $ called
the Planck constant. For each $\hbar $ the protocol associates to
observables from $\mathcal{A}$ and automorphisms from $\Gamma $ certain
operators acting on a finite dimensional Hilbert space $\mathcal{H}$. \ 

Let $\hbar =\frac{1}{p},$ where $p$ is an odd prime number, and consider the
additive character $\psi :\mathbb{F}_{p}\longrightarrow {\mathbb{C}}^{\ast
},\; \psi (t)=e^{\frac{2\pi i}{p}t}$. Define the non-commutative torus \cite%
{Co, Ri} $\mathcal{A}_{\hbar }$ to be the free non-commutative $%
%TCIMACRO{\U{2102} }%
%BeginExpansion
\mathbb{C}
%EndExpansion
$ -algebra generated by the symbols $s(\xi )$, $\xi \in \Lambda ,$ and the
relations 
\begin{equation}
\psi (\tfrac{1}{2}\omega (\xi ,\eta ))s(\xi +\eta )=s(\xi )s(\eta ).
\label{CCR}
\end{equation}%
Here we consider $\omega $ as a map $\omega :\Lambda \times \Lambda
\longrightarrow \mathbb{F}_{p}$.

Note that $\mathcal{A}_{\hbar }$ satisfies the following properties:

\begin{itemize}
\item As a vector space $\mathcal{A}_{\hbar }$ is equipped with a natural
basis $\ s(\xi )$, $\xi \in \Lambda $. \ Hence we can identify the vector
space $\mathcal{A}_{\hbar }$ with the vector space $\mathcal{A}$ for each
value of $\hbar ,$ 
\begin{equation}
\mathcal{A\simeq A}_{\hbar }.  \label{ident}
\end{equation}

\item Substituting $\hbar =0$ we have $\mathcal{A}=\mathcal{A}_{0}$. Hence,
we see that indeed $\mathcal{A}_{\hbar }$ is a deformation of the algebra of
trigonometric polynomials on $\mathbb{T}$.

\item The group $\Gamma $ acts by automorphisms on the algebra $\mathcal{A}%
_{\hbar }$, via $\gamma \cdot s(f)=s(\gamma f)$, where $\gamma \in \Gamma $
and $f\in \mathcal{A}_{\hbar }$. This action induces an action of $\Gamma $
on the category of representations of $\mathcal{A}_{\hbar }$, taking a
representation $\pi $ and sending it to the representation $\pi ^{\gamma }$,
where $\pi ^{\gamma }(f)=\pi (\gamma f)$.
\end{itemize}

We use the identification $\left( \text{\ref{ident}}\right) $ and a
distinguished representation of the algebra $\mathcal{A}_{\hbar }$ to
describe the quantization of the functions. All the irreducible algebraic
representations of $\mathcal{A}_{\hbar }$ are classified \cite{GH2} and each
of them is of dimension $p^{N}.$

\begin{theorem}[Invariant representation \protect \cite{GH2}]
\label{fixrep}Let $\hbar =\frac{1}{p}$ where $p$ is a prime number. There
exists a unique (up to isomorphism) irreducible representation $\pi :%
\mathcal{A}_{\hbar }\rightarrow End(\mathcal{H}_{\hbar })$ which is fixed by
the action of $\Gamma $. Namely, $\pi ^{\gamma }$ is isomorphic to $\pi $
for every $\gamma \in \Gamma $.
\end{theorem}

Let $(\pi ,\mathcal{A}_{\hbar },\mathcal{H)}$ be a representative of the the
special representation defined in Theorem \ref{fixrep}. For every element $%
\gamma \in \Gamma $ we have an isomorphism $\widetilde{\rho }(\gamma ):%
\mathcal{H\rightarrow H}$ intertwining the representations $\pi $ and $\pi
^{\gamma }$, namely, it satisfies 
\begin{equation}
\widetilde{\rho }(\gamma )\pi (f)\widetilde{\rho }(\gamma )^{-1}=\pi (\gamma
f),  \label{Ego}
\end{equation}%
for every $f\in \mathcal{A}_{\hbar }$ and $\gamma \in \Gamma $. The
isomorphism $\widetilde{\rho }(\gamma )$ is not unique but unique up to a
scalar (this is a consequence of Schur's lemma). It is easy to realize that
the collection $\{ \widetilde{\rho }(\gamma )\}$ constitutes a projective
representation $\widetilde{\rho }:\Gamma \rightarrow PGL(\mathcal{H)}$.
Assume now that $\hbar =\frac{1}{p}$ where $p$ is an odd prime $\neq 3.$We
have the following linearization theorem:

\begin{theorem}[Linearization \protect \cite{GH3, GH9}]
\label{factorization}There exist a unique representation $\rho :\Gamma
\rightarrow GL(\mathcal{H)}$ that that satisfies (\ref{Ego}) and factors
through the quotient group $Sp\simeq Sp(2N,\mathbb{F}_{p})$.
\end{theorem}

Concluding, we established a distinguished pair of representations%
\begin{equation*}
\rho :\Gamma \rightarrow GL(\mathcal{H})\text{ \  \ and \  \ }\pi :\mathcal{A}%
_{\hbar }\rightarrow End(\mathcal{H}),
\end{equation*}%
satisfying the compatibility condition (\ref{Ego}).

\begin{remark}
The representation $\rho :Sp$ $\rightarrow GL(\mathcal{H)}$ defined by
Theorem \ref{factorization} is the Weil representation associated with the
additive character $\psi (t)=e^{\frac{2\pi i}{p}t}$---here obtained in a
different manner via quantization of the torus.
\end{remark}

\subsection{ The quantum dynamical system}

Recall that we started with a dynamic on $\mathbb{T}$, generated by a
generic (i.e., regular with no non-trivial invariant co-isotropic sub-tori)
element $A\in \Gamma $. Using the Weil representation we can associate to $A$
the unitary operator $\rho (A):\mathcal{H\rightarrow H},$ which constitutes
the generator of discrete time quantum dynamics. \ We would like to study $%
\rho (A)$-eigenstates$\ $%
\begin{equation}
\rho (A)\varphi =\lambda \varphi ,\text{ \  \ }\varphi \in \mathcal{H}\text{,}
\label{EP}
\end{equation}%
which satisfy addition arithmetic symmetries.

\section{Hecke quantum unique ergodicity\label{HQUE}}

It turns out that the operator $\rho (A)$ has degeneracies---its eigenspaces
might be extremely large. This is manifested in the existence of a group of
hidden\footnote{%
The group $Z(A,\Gamma ),$ of \ the elements in $\Gamma $ that commute with $%
A,$ does not contribute much to the harmonic analysis of $\rho (A).$}
symmetries commuting with $\rho (A).$ These symmetries can be computed using
the Weil representation. Indeed, let $T_{A}=Z(A,Sp)$ be the centralizer of
the element $A$ in the group $Sp$. Clearly $T_{A}$ contains the cyclic group 
$\left \langle A\right \rangle $ generated by the element $A,$ but it often
happens that $T_{A}$ contains additional elements. The assumption that $A$
is generic implies that for sufficiently large $p$ (so that $p$ does not
divides the discriminant of $A$) the group $T_{A}$ consists of the $\mathbb{F%
}_{p}$-rational points of a maximal torus $\mathbf{T}_{A}\subset \mathbf{Sp,}
$ i.e., $T_{A}=\mathbf{T}_{A}(\mathbb{F}_{p}).$ We will call the group $%
T_{A} $ the \textit{Hecke }torus. It acts semisimply on $\mathcal{H}$,
decomposing it into a direct sum of character spaces $\mathcal{H=}\oplus 
\mathcal{H}_{\chi }$ where $\chi $ runs in the group of characters of $T_{A}$%
. We shall study common eigenstates $\varphi \in $ $\mathcal{H}_{\chi },$
which we will call in this setting \textit{Hecke eigenstates }and will be
assumed to be normalized so that $\left \Vert \varphi \right \Vert _{%
\mathcal{H}}=1$. In particular, we will bound the Wigner distributions $%
\left \langle \varphi |\pi (f)\varphi \right \rangle ,$ where $f$ $\in 
\mathcal{A}$ is an observable on the torus $\mathbb{T}$. We will call these
matrix coefficients \textit{Hecke--Wigner distributions}. It will be
convenient for us to treat two cases.

\subsection{The strongly generic case\label{TSGC}}

Let us assume first that the automorphism $A$ acts on $\mathbb{T}$ with no
invariant sub-tori. In dual terms, this means that the element $A$ acts
irreducibly on the $%
%TCIMACRO{\U{211a} }%
%BeginExpansion
\mathbb{Q}
%EndExpansion
$-vector space $\Lambda _{%
%TCIMACRO{\U{211a} }%
%BeginExpansion
\mathbb{Q}
%EndExpansion
}=\Lambda \otimes _{%
%TCIMACRO{\U{2124} }%
%BeginExpansion
\mathbb{Z}
%EndExpansion
}%
%TCIMACRO{\U{211a} }%
%BeginExpansion
\mathbb{Q}
%EndExpansion
.$ We denote by $r_{p}$ the symplectic rank of $T_{A}$, \ i.e., $r_{p}=|\Xi
| $ where $\Xi =\Xi (T_{A})$ is the symplectic type of $T_{A}$ (see
Definition \ref{type}). By definition we have $1\leq r_{p}\leq N$.

\begin{theorem}
\label{SGC}Consider a non-trivial exponent $0\neq \xi \in \Lambda $ and a
sufficiently large prime number $p.$ Then for every unit Hecke eigenstate $%
\varphi \in $ $\mathcal{H}_{\chi }$ the following bound holds: 
\begin{equation}
\left \vert \left \langle \varphi |\pi (\xi )\varphi \right \rangle \right
\vert \leq \frac{m_{\chi }\cdot (2+o(1))^{r_{p}}}{\sqrt{p}^{N}},
\label{SGCF}
\end{equation}%
where $m_{\chi }=\dim \mathcal{H}_{\chi }.$
\end{theorem}

The lattice $\Lambda $ constitutes a basis for $\mathcal{A};$ hence, using
the bound (\ref{SGCF}) we obtain:

\begin{theorem}[Hecke quantum unique ergodicity---strongly generic case]
Consider an observable $f\in \mathcal{A}$ and a sufficiently large prime
number $p.$ Then for every normalized Hecke eigenstate $\varphi $ we have 
\begin{equation*}
\left \vert \left \langle \varphi |\pi (f)\varphi \right \rangle -\tint
\limits_{\mathbb{T}}fd\mu \right \vert \leq \frac{C_{f}}{\sqrt{p}^{N}},
\end{equation*}%
where $\mu =|\omega |^{N}$ is the corresponding volume form and $C_{f}$ is
an explicit computable constant which depends only on the function $f.$
\end{theorem}

\begin{remark}
In Subsection \ref{St} we will elaborate on the distribution of the
symplectic rank $r_{p}$ $\left( \text{\ref{SGCF}}\right) $.
\end{remark}

\subsubsection{Proof of Theorem \protect \ref{SGC}}

The proof is by reduction to the bound on the Hecke--Wigner distributions
obtained in Section \ref{B}, i.e., reduction to Theorem \ref{KR_general1}.
Our first goal is to interpret the Hecke--Wigner distribution $\left \langle
\varphi |\pi (\xi )\varphi \right \rangle $ in terms of the Heisenberg--Weil
representation.

\textbf{Step 1.} \textit{Replacing the non-commutative torus by the finite
Heisenberg group.} Note that the Hilbert space $\mathcal{H}$ is a
representation space of both the algebra $\mathcal{A}_{\hbar }$ and the
group $Sp$ via $\pi $ and $\rho ,$ respectively. We will show next that the
representation $(\pi ,\mathcal{A}_{\hbar },\mathcal{H})$ is
\textquotedblleft equivalent\textquotedblright \ to the Heisenberg
representation of some finite Heisenberg group. The representation $\pi $ is
determined by its restriction to the lattice $\Lambda $. However, the
restriction $\pi _{|\Lambda }:\Lambda \rightarrow GL(\mathcal{H)}$ is not an
homomorphism and in fact constitutes (see formula (\ref{CCR})) a projective
representation of the lattice given by $\psi (\tfrac{1}{2}\omega (\xi ,\eta
))\pi (\xi +\eta )=\pi (\xi )\pi (\eta )$. It is evident from this formula
that the map $\pi _{|\Lambda }$ factors through the quotient $\mathbb{F}_{p}$%
-vector space $V=\Lambda /p\Lambda ,$ i.e., $\Lambda \rightarrow V=\Lambda
/p\Lambda \rightarrow GL(\mathcal{H)}.$ The vector space $V$ is equipped
with a symplectic structure $\omega $ obtained via specialization of the
form on $\Lambda .$ Let $H=H(V)$ be the Heisenberg group associated with $%
(V,\omega )$. So the map $\pi :V\rightarrow GL(\mathcal{H)}$ lifts into an
honest representation of the Heisenberg group $\pi :H\rightarrow GL(\mathcal{%
H)}.$ Finally, the Heisenberg representation $\pi $ and the Weil
representation $\rho $ glue into a single representation $\tau =\rho \ltimes
\pi $ of the Jacobi group $J=Sp\ltimes H$, which is of course nothing other
than the Heisenberg--Weil representation 
\begin{equation}
\tau :J\rightarrow GL(\mathcal{H)}.  \label{HWR}
\end{equation}

\textbf{Step 2. }\textit{Reformulation. }Let $\mathbf{V}$ and $\mathbf{T}%
_{A} $ be the algebraic group scheme defined over $%
%TCIMACRO{\U{2124} }%
%BeginExpansion
\mathbb{Z}
%EndExpansion
$ so that $\Lambda =\mathbf{V(%
%TCIMACRO{\U{2124} }%
%BeginExpansion
\mathbb{Z}
%EndExpansion
)}$ and for every prime $p$ we have $V=\mathbf{V}(\mathbb{F}_{p})$ and $%
T_{A} $ $=\mathbf{T}_{A}(\mathbb{F}_{p}).$ In this setting for every prime
number $p$ we can consider the lattice element $\xi \in \Lambda $ as a
vector in the $\mathbb{F}_{p}$-vector space $V$.

Let $(\tau ,J,\mathcal{H})$ be the Heisenberg--Weil representation (\ref{HWR}%
) and consider a unit Hecke eigenstate $\varphi \in \mathcal{H}_{\chi }.$ We
need to verify that for a sufficiently large prime number $p$ we have 
\begin{equation}
\left \vert \left \langle \varphi |\pi (\xi )\varphi \right \rangle \right
\vert \leq \frac{m_{\chi }\cdot (2+o(1))^{r_{p}}}{\sqrt{p}^{N}},  \label{ref}
\end{equation}%
where $m_{\chi }=\dim \mathcal{H}_{\chi }$ and $r_{p}$ is the symplectic
rank of $T_{A}.$ \smallskip

\textbf{Step 3.} \textit{Verification. }We need to show that we meet the
conditions of Theorem \ref{KR_general1}. What is left to check is that for
sufficiently large prime number $p$ the vector $\xi \in V$ is not contained
in any $T_{A}$-invariant subspace of $V.$ Let us denote by $O_{\xi }$ the
orbit $O_{\xi }=T_{A}\cdot \xi .$ We need to show that for a sufficiently
large $p$ we have 
\begin{equation}
Span_{\mathbb{F}_{p}}\{O_{\xi }\}=V.  \label{Ir}
\end{equation}

The condition (\ref{Ir}) is satisfied since it holds globally. In more
details, our assumption on $A$ guarantees that it holds for the
corresponding objects over the field of rational numbers $%
%TCIMACRO{\U{211a} }%
%BeginExpansion
\mathbb{Q}
%EndExpansion
$, i.e., $Span_{%
%TCIMACRO{\U{211a} }%
%BeginExpansion
\mathbb{Q}
%EndExpansion
}\{ \mathbf{T}_{A}(%
%TCIMACRO{\U{211a} }%
%BeginExpansion
\mathbb{Q}
%EndExpansion
)\cdot \xi \}=\mathbf{V}\left( 
%TCIMACRO{\U{211a} }%
%BeginExpansion
\mathbb{Q}
%EndExpansion
\right) .$ Hence (\ref{Ir}) holds for sufficiently large prime number $p.$

This completes the proof of Theorem \ref{SGC}.

\subsection{The distribution of the symplectic rank\label{St}}

We would like to compute the asymptotic distribution of the symplectic rank $%
r_{p}$ (\ref{ref}) in the set $\left \{ 1,...,N\right \} $, \ i.e., 
\begin{equation}
\delta (r)=\lim_{x\rightarrow \infty }\frac{\# \left \{ r_{p}=r\text{ };%
\text{ }p\leq x\right \} }{\pi (x)},  \label{delta}
\end{equation}%
where $\pi (x)$ denotes the number of prime numbers up to $x.$

We fix an algebraic closure $\overline{%
%TCIMACRO{\U{211a} }%
%BeginExpansion
\mathbb{Q}
%EndExpansion
}$ of the field $%
%TCIMACRO{\U{211a} }%
%BeginExpansion
\mathbb{Q}
%EndExpansion
,$ and denote by $G$ the Galois group $G=Gal(\overline{%
%TCIMACRO{\U{211a} }%
%BeginExpansion
\mathbb{Q}
%EndExpansion
}/%
%TCIMACRO{\U{211a} }%
%BeginExpansion
\mathbb{Q}
%EndExpansion
).$ Consider the vector space $\mathbf{V}=\mathbf{V(}\overline{%
%TCIMACRO{\U{211a} }%
%BeginExpansion
\mathbb{Q}
%EndExpansion
})$. By extension of scalars the symplectic form $\omega $ on $\mathbf{V}%
\left( 
%TCIMACRO{\U{211a} }%
%BeginExpansion
\mathbb{Q}
%EndExpansion
\right) $ induces a $\overline{%
%TCIMACRO{\U{211a} }%
%BeginExpansion
\mathbb{Q}
%EndExpansion
}$-linear symplectic form on $\mathbf{V}$ which we will also denote by $%
\omega .$ Let $\mathbf{T}$ denote the algebraic torus $\mathbf{T=}$ $\mathbf{%
T}_{A}(\overline{%
%TCIMACRO{\U{211a} }%
%BeginExpansion
\mathbb{Q}
%EndExpansion
}).$ The action of $\mathbf{T}$ on $\mathbf{V}$ is completely reducible,
decomposing it into one-dimensional character spaces $\mathbf{V=}\oplus
_{\chi \in \mathfrak{X}}\mathbf{V}_{\chi }.$ Let $\Theta $ be the
restriction of the symplectic transpose $(\cdot )^{t}:End(\mathbf{V}%
)\rightarrow End(\mathbf{V})$ to $\mathbf{T.}$ The involution $\Theta $ acts
on the set of characters $\mathfrak{X}$ by $\chi \mapsto \Theta (\chi )=\chi
^{-1}$ and this action is compatible with the action of the Galois group $G$
on $\mathfrak{X}$ by conjugation $\chi \mapsto g\chi g^{-1}$, where $\chi
\in \mathfrak{X}$ and $g\in G$. This means (recall that $A$ is strongly
generic) that we have a transitive action of $G$ on the set $\mathfrak{X}%
/\Theta .$ Consider the kernel $\mathrm{K}=\ker (G\rightarrow Aut(\mathfrak{X%
}/\Theta ))$ and the corresponding finite Galois group $Q=G/\mathrm{K}.$
Considering $Q$ as a subgroup of $\ Aut(\mathfrak{X}/\Theta )$ we define the
cycle number $c(C)$ of a conjugacy class $C\subset Q$ to be the number of
irreducible cycles that compose a representative of $C.$ A direct
application of the Chebotarev density theorem \cite{Se} is the following:

\begin{proposition}[Chebotarev's theorem]
\label{CT}The distribution $\delta $ $\left( \text{\ref{delta}}\right) $
obeys%
\begin{equation*}
\delta (r)=\frac{\left \vert C_{r}\right \vert }{\left \vert Q\right \vert },
\end{equation*}%
where $C_{r}=\underset{c(C)=r}{\underset{C\subset Q}{\cup }}C.$
\end{proposition}

For a proof see Appendix \ref{PCT}.

\subsection{The general generic case \label{GGC}}

Let us now treat the more general case where the automorphism $A$ acts on $%
\mathbb{T}$ in a generic way (Definition \  \ref{qergodic}). In dual terms,
this means that the torus $\mathbf{T}(%
%TCIMACRO{\U{211a} }%
%BeginExpansion
\mathbb{Q}
%EndExpansion
)=\mathbf{T}_{A}(%
%TCIMACRO{\U{211a} }%
%BeginExpansion
\mathbb{Q}
%EndExpansion
)$ acts on the symplectic vector space $\mathbf{V}(%
%TCIMACRO{\U{211a} }%
%BeginExpansion
\mathbb{Q}
%EndExpansion
)$ $=\Lambda \otimes _{%
%TCIMACRO{\U{2124} }%
%BeginExpansion
\mathbb{Z}
%EndExpansion
}%
%TCIMACRO{\U{211a} }%
%BeginExpansion
\mathbb{Q}
%EndExpansion
$ decomposing it into an orthogonal symplectic direct sum%
\begin{equation}
(\mathbf{V}(%
%TCIMACRO{\U{211a} }%
%BeginExpansion
\mathbb{Q}
%EndExpansion
),\omega )=\tbigoplus \limits_{\alpha \in \Xi }(\mathbf{V}_{\alpha }(%
%TCIMACRO{\U{211a} }%
%BeginExpansion
\mathbb{Q}
%EndExpansion
),\omega _{\alpha }),  \label{Qdecomp}
\end{equation}%
with an irreducible action of $\mathbf{T}(%
%TCIMACRO{\U{211a} }%
%BeginExpansion
\mathbb{Q}
%EndExpansion
)$ on each of the spaces $\mathbf{V}_{\alpha }(%
%TCIMACRO{\U{211a} }%
%BeginExpansion
\mathbb{Q}
%EndExpansion
)$. For an element $\xi \in \Lambda $ define its support with respect to the
decomposition (\ref{Qdecomp}) by $S_{\xi }=Supp(\xi )=\{ \alpha ;$ $%
P_{\alpha }\xi \neq 0\},$ where $P_{\alpha }:\mathbf{V}(%
%TCIMACRO{\U{211a} }%
%BeginExpansion
\mathbb{Q}
%EndExpansion
)\rightarrow $ $\mathbf{V}(%
%TCIMACRO{\U{211a} }%
%BeginExpansion
\mathbb{Q}
%EndExpansion
)$ is the projector onto the space $\mathbf{V}_{\alpha }(%
%TCIMACRO{\U{211a} }%
%BeginExpansion
\mathbb{Q}
%EndExpansion
)$ and denote by $d_{\xi }$ the dimension $d_{\xi }=\sum_{\alpha \in S_{\xi
}}\dim \mathbf{V}_{\alpha }(%
%TCIMACRO{\U{211a} }%
%BeginExpansion
\mathbb{Q}
%EndExpansion
).$ The decomposition (\ref{Qdecomp}) induces a decomposition of the torus $%
\mathbf{T}(%
%TCIMACRO{\U{211a} }%
%BeginExpansion
\mathbb{Q}
%EndExpansion
)$ into a product of completely inert tori 
\begin{equation}
\mathbf{T}(%
%TCIMACRO{\U{211a} }%
%BeginExpansion
\mathbb{Q}
%EndExpansion
)=\tprod \limits_{\alpha \in \Xi }\mathbf{T}_{\alpha }(%
%TCIMACRO{\U{211a} }%
%BeginExpansion
\mathbb{Q}
%EndExpansion
).  \label{TDQ}
\end{equation}%
Considering now a sufficiently large prime number $p$ and specialize all the
object involved to the finite field $\mathbb{F}_{p}$. The Hecke torus $T=%
\mathbf{T}(\mathbb{F}_{p})$ acts via the Weil representation on the Hilbert
space $\mathcal{H}\mathbb{\ }$decomposing it into an orthogonal direct sum $%
\mathcal{H=}\oplus \mathcal{H}_{\chi }.$ The decomposition (\ref{TDQ})
induces decompositions on the level of groups of points $T=\tprod T_{\alpha
},$ where $T_{\alpha }=\mathbf{T}_{\alpha }(\mathbb{F}_{p}),$ on the level
of characters $\chi =\Pi \chi _{\alpha }:\tprod T_{\alpha }\rightarrow 
%TCIMACRO{\U{2102} }%
%BeginExpansion
\mathbb{C}
%EndExpansion
^{\ast }$, and on the level of character spaces $\mathcal{H}_{\chi }=\oplus 
\mathcal{H}_{\chi _{\alpha }}.$ For each torus $T_{\alpha }$ we denote by $%
r_{p,\alpha }=r_{p}(T_{\alpha })$ its symplectic rank (see Definition \ref%
{type}) and we consider the integer $\left \vert S_{\xi }\right \vert \leq
r_{p,\xi }\leq d_{\xi }$ given by $r_{p,\xi }=\underset{\alpha \in S_{\xi }}{%
\Pi }r_{p,\alpha }.$ Let us denote by $m_{\chi _{\xi }}$ the dimension $%
m_{\chi _{\xi }}=\sum_{\alpha \in S_{\xi }}\dim \mathcal{H}_{\chi _{\alpha
}}.$

\begin{theorem}
\label{GCT}Consider a non-trivial exponent $0\neq \xi \in \Lambda $ and a
sufficiently large prime number $p.$ Then for every unit Hecke eigenstate $%
\varphi \in $ $\mathcal{H}_{\chi }$ we have 
\begin{equation*}
\left \vert \left \langle \varphi |\pi (\xi )\varphi \right \rangle \right
\vert \leq \frac{m_{\chi _{\xi }}\cdot (2+o(1))^{r_{p,\xi }}}{\sqrt{p}%
^{d_{\xi }}}.
\end{equation*}
\end{theorem}

We consider the decomposition (\ref{Qdecomp}) and denote by $d=\min_{\alpha }%
\mathbf{V}_{\alpha }(%
%TCIMACRO{\U{211a} }%
%BeginExpansion
\mathbb{Q}
%EndExpansion
).$ Using the lattice $\Lambda $ as a basis for the algebra $\mathcal{A}$ we
obtain:

\begin{theorem}[Hecke quantum unique ergodicity---generic case]
Consider an observable $f\in \mathcal{A}$ and a sufficiently large prime
number $p.$ Then for every unit Hecke eigenstate $\varphi $ we have 
\begin{equation*}
\left \vert \left \langle \varphi |\pi (f)\varphi \right \rangle -\tint
\limits_{\mathbb{T}}fd\mu \right \vert \leq \frac{C_{f}}{\sqrt{p}^{d}},
\end{equation*}%
where $\mu =|\omega |^{N}$ is the corresponding volume form and $C_{f}$ is
an explicit computable constant which depends only on the function $f.$
\end{theorem}

The proof of Theorem \ref{GCT} is a straightforward application of Theorem %
\ref{SGC}. In more details, considering the \textquotedblleft
global\textquotedblright \ decomposition (\ref{Qdecomp}) of the torus $%
\mathbf{T}(%
%TCIMACRO{\U{211a} }%
%BeginExpansion
\mathbb{Q}
%EndExpansion
)$ to a product of completely inert tori $\mathbf{T}_{\alpha }(%
%TCIMACRO{\U{211a} }%
%BeginExpansion
\mathbb{Q}
%EndExpansion
)$ we may apply the theory developed for the strongly generic case in
Subsection (\ref{TSGC}) to each of the tori $\mathbf{T}_{\alpha }(%
%TCIMACRO{\U{211a} }%
%BeginExpansion
\mathbb{Q}
%EndExpansion
)$ to deduce Theorem \ref{GCT}.

\begin{remark}
As explained in Subsection $\left( \text{\ref{St}}\right) $ the distribution
of the symplectic rank $r_{p,\xi }$ is determined by the Chebotarev theorem
applied to $\left( \text{now a product of}\right) $ suitable finite Galois
groups $Q_{\alpha \text{ }}$attached to the tori $\mathbf{T}_{\alpha },$ $%
\alpha \in S_{\xi }$ $\left( \text{\ref{TDQ}}\right) .$
\end{remark}

\begin{remark}
The corresponding quantum unique ergodicity theorem for statistical states
of generic automorphism $A$ of $\mathbb{T}$ $($cf. Theorem \ref{QUESS}$)$
follows directly from Theorem \ref{GCT} .
\end{remark}

\appendix

\section{Proofs\label{Pr}}

\subsection{Proof of Theorem \protect \ref{Multiplicities}\label{PM}}

Fix a character $\chi $ of $T$ and denote by $P_{\chi }$ the orthogonal
projector on the character space $\mathcal{H}_{\chi }.$ The projector $%
P_{\chi }$ can be written in terms of the Weil representation $P_{\chi
}=|T|^{-1}\sum_{g\in T}\chi ^{-1}(g)\rho (g);$ therefore $m_{\chi }=\dim 
\mathcal{H}_{\chi }=Tr(P_{\chi })=|T|^{-1}\sum_{g\in T}\chi ^{-1}(g)Tr(\rho
(g))$ which by the character formula (\ref{Char W}) is equal to 
\begin{equation*}
\frac{1}{|T|}\left( \left( \dsum \limits_{g\in T\smallsetminus I}\chi
^{-1}(g)\sigma (-\det (g-I))\right) +q\right) .
\end{equation*}%
From the orthogonality relations for characters, and the fact that $|T|=q-1$
if $T$ splits and equal $q+1$ if $T$ is inert, we see that it is enough to
prove the following claim:

\begin{claim}
\label{Rest}On $T\smallsetminus I$ \ we have $\sigma (-\det (g-I))=\pm
\sigma _{T}(g)$ where the $+$ and $-$ are attained when the torus $T$ is
split or inert, respectively. \smallskip
\end{claim}

\subsubsection{Proof of Claim \protect \ref{Rest}}

The equality is easily satisfied in the split case. Let us assume that $T$
is inert and identify $Sp\simeq SL(2,\mathbb{F}_{q})\subset SL(2,\mathbb{F}%
_{q^{2}})$. There exists a matrix $S\in $ $SL(2,\mathbb{F}_{q^{2}})$ so that 
\begin{equation*}
STS^{-1}=\left \{ \left( 
\begin{array}{cc}
c & 0 \\ 
0 & c^{-1}%
\end{array}%
\right) ;\text{ }c\in C\right \} ,
\end{equation*}%
where $C$ is the kernel of the norm map $C=Ker(N:\mathbb{F}%
_{q^{2}}\rightarrow \mathbb{F}_{q})$ $\subset \mathbb{F}_{q^{2}}.$ In these
coordinates our claim reduced to the identity $(\frac{(c-1)^{2}}{c})^{\frac{%
q-1}{2}}=-c^{\frac{q+1}{2}}.$ This completes the proof of the claim and of
Theorem \ref{Multiplicities}.

\subsection{Proof of Proposition \protect \ref{prop_dim}, and Claim \protect
\ref{Comm}\label{ComDim}}

Consider a symplectic vector space $(V,\omega )$ over $k$ and a maximal
torus $T\subset Sp(V,\omega ).$ Assume that $T$ acts irreducibly on $V.$
Denote by $A=Z(T,End(V))$ the centralizer of $T$ in the algebra of all
linear endomorphisms of $V.$ We need to show that $A$ is commutative and
that $\dim _{K}V=2$, where $K=A^{\Theta },$ with $\Theta $ being the
canonical quadratic element in $Gal(A/k)$ obtained by restricting the
symplectic transpose defined by $\omega $ on $End(V)$ to $A.$

Let $\overline{k}$ denote an algebraic closure of the field $k$ and denote
by $G$ the Galois group $G=Gal(\overline{k}/k)$. Consider the vector space $%
\mathbf{V}=\overline{k}\otimes V$. By extension of scalars the symplectic
structure $\omega $ induces a $\overline{k}$-linear symplectic structure on $%
\mathbf{V}$, which we denote also by $\omega .$ Let $\mathbf{T}$ denote the
algebraic torus, i.e., $\mathbf{T=}\overline{k}\otimes T$ . Consider the
algebra $\mathbf{A=}Z(\mathbf{T,}End(\mathbf{V}))$. Note that in this
setting $\mathbf{A}$ is not necessarily a division algebra. Let $\Theta $ be
the restriction of the symplectic transpose $(\cdot )^{t}:End(\mathbf{V}%
)\rightarrow End(\mathbf{V}),$ to the algebra $\mathbf{A}$ and denote by $%
\mathbf{K=A}^{\Theta }$ the subalgebra consisting of elements $a\in \mathbf{A%
}$ which are fixed by $\Theta .$ The group $G$ acts on all structures
involved. We have $V=\mathbf{V}^{G},$ $T=\mathbf{T}^{G},$ $A=\mathbf{A}^{G}$
and $K=\mathbf{K}^{G}.$

The action of $\mathbf{T}$ on $\mathbf{V}$ is completely reducible,
decomposing it into one-dimensional character spaces 
\begin{equation}
\mathbf{V=}\tbigoplus \limits_{\chi \in \mathfrak{X}}\mathbf{V}_{\chi }.
\label{decomposition1}
\end{equation}

The set $\mathfrak{X}$, consists of $N$ pairs of characters, $\chi ,\chi
^{-1}\in \mathfrak{X}$. The algebra $\mathbf{A}$ consists of operators which
are diagonal with respect to the decomposition (\ref{decomposition1}) $%
\mathbf{A}=diag(a_{\chi }\in \overline{k}:\chi \in \mathfrak{X}).$ In
particular, this implies that $A=\mathbf{A}^{G}$ is commutative as was
claimed. The involution $\Theta $ can be described as $a=diag(a_{\chi
})\longmapsto \Theta (a)=diag(\Theta (a)_{\chi }=a_{\chi ^{-1}}).$
Therefore, we obtain that $\mathbf{K}=\left \{ a=diag(a_{\chi }):a_{\chi
}=a_{\chi ^{-1}}\right \} .$

Finally, we observe that $\mathbf{V}$ is a free module of rank $2$ over the
algebra $\mathbf{K}$. By Hilbert's Theorem 90, this property continues to
hold after taking Galois invariants, that is we obtain that $\dim _{K}V=2$.

\subsection{Proof of Corollary \protect \ref{DT}\label{Norm}}

We need to show that $T=\left \{ a\in A:N_{A/K}(a)=1\right \} .$ Because $T$
is a maximal torus in $Sp$, it coincides with its own centralizer $Z(T,Sp)$.
We have $Z(T,Sp)=A\cap Sp$. An element $a\in A$ lies in the group $Sp$ if it
satisfies $\omega (au,av)=\omega (u,v)$ for every $u,v\in V$. This is
equivalent to $\omega (u,\Theta (a)av)=\omega (u,v)$ which in turn implies
that $N_{A/K}(a)=\Theta (a)a=1$.

\subsection{Proof of Proposition \protect \ref{canonicallift}\label{C-L}}

We should prove that there exists a $T$-invariant $K$-linear symplectic form 
$\overline{\omega }:V\times V\rightarrow K$, satisfying the property $%
Tr\circ \overline{\omega }=\omega .$ Consider the decomposition (\ref%
{decomposition1}) $\mathbf{V=}\oplus \mathbf{V}_{\chi }.$ Define $\overline{%
\omega }:\mathbf{V\times V\rightarrow K}$ as follows 
\begin{equation*}
\overline{\omega }(\sum v_{\chi },\sum u_{\chi })=diag(a_{\chi }=\omega
(v_{\chi },u_{\chi ^{-1}})+\omega (v_{\chi ^{-1}},u_{\chi })).
\end{equation*}%
We have $Tr\circ \overline{\omega }=\omega .$ Clearly, the form $\overline{%
\omega }$ is invariant under the action of the torus $\mathbf{T}$. \
Finally, the form $\overline{\omega }$ commutes with the Galois action,
hence it restricts to give a desired form $\overline{\omega }:V\times
V\rightarrow K$, \ which is $T$-invariant and satisfies $Tr\circ \overline{%
\omega }=\omega $.

\subsection{Proof of Lemma \protect \ref{technical} \label{Ptechnical}}

We use the notation of Subsection \ref{ComDim}. Consider the decomposition $%
\mathbf{V=}\oplus _{\chi \in \mathfrak{X}}\mathbf{V}_{\chi }.$ The Galois
group $G=Gal(\overline{k}/k)$ acts on the set of characters $\mathfrak{X}$.
The action is given by conjugation $\chi \mapsto g\chi g^{-1}$, where $\chi
\in \mathfrak{X}$ and $g\in G$. The set $\mathfrak{X}$ decomposes into a
union of $G$ orbits, $\mathfrak{X=}\tbigcup \mathcal{O}_{\beta }$. To each
orbit $\mathcal{O}_{\beta }$ there exists a unique dual orbit $\mathcal{O}_{%
\widehat{\beta }}$ so that $\chi \in $ $\mathcal{O}_{\beta }$ if and only if 
$\chi ^{-1}\in \mathcal{O}_{\widehat{\beta }}$ (note that sometimes $%
\mathcal{O}_{\widehat{\beta }}=\mathcal{O}_{\beta }$). Let $\  \mathcal{O}%
_{\alpha }$ denote the union $\mathcal{O}_{\alpha }=\mathcal{O}_{\beta }\cup 
\mathcal{O}_{\widehat{\beta }}$. We denote by $\Xi $ the set of $\mathcal{O}%
_{\alpha }$'s. We use the following terminology: If $\mathcal{O}_{\widehat{%
\beta }}=\mathcal{O}_{\beta }$,\ we say that $\mathcal{O}_{\alpha }$ is of
type I, otherwise we say that $\mathcal{O}_{\alpha }$ is of type II.

The decomposition $\mathfrak{X=}\tbigcup \limits_{\alpha \in \Xi }\mathcal{O}%
_{\alpha }$ induces a decomposition 
\begin{equation}
(\mathbf{V,\omega )=}\tbigoplus \limits_{\alpha \in \Xi }(\mathbf{V}_{\alpha
},\omega _{\alpha }),  \label{decom1}
\end{equation}%
where $\mathbf{V}_{\alpha }=\oplus _{\chi \in \mathcal{O}_{\alpha }}\mathbf{V%
}_{\chi }$. This in turn induces the following decompositions: 
\begin{equation}
\mathbf{T}\mathbf{=}\tprod \mathbf{T}_{\alpha },\text{ }\mathbf{A}\mathbf{=}%
\tbigoplus \mathbf{A}_{\alpha },\text{\ }  \label{dec2}
\end{equation}%
where $\mathbf{T}_{\alpha }$ is the sub-torus consisting of elements $g\in 
\mathbf{T}$ such that $g_{|\mathbf{V}_{\beta }}=Id$ for every $\beta \neq
\alpha $. In the same spirit $\mathbf{A}_{\alpha }$ is the subalgebra
consisting of elements $a\in \mathbf{A}$ such that $a_{|\mathbf{V}_{\beta
}}=Id$, for every $\beta \neq \alpha $. It is easy to verify that each
algebra $\mathbf{A}_{\alpha }$ is closed under the involution $\Theta $. Let 
$\mathbf{K}_{\alpha }$ denote the invariant subalgebra $\mathbf{K}_{\alpha }$
$=\mathbf{A}_{\alpha }^{\Theta }$. We have the following decomposition: 
\begin{equation}
\mathbf{K=}\tbigoplus \mathbf{K}_{\alpha }.  \label{decom3}
\end{equation}

Finally, following the same construction as in Subsection \ref{C-L}, we can
lift the symplectic form $\omega _{\alpha }$ to a $\mathbf{K}_{\alpha }$%
-linear $\mathbf{T}_{\alpha }$-invariant symplectic form 
\begin{equation}
\overline{\omega }_{\alpha }:\mathbf{V}_{\alpha }\times \mathbf{V}_{\alpha
}\rightarrow \mathbf{K}_{\alpha },  \label{decom4}
\end{equation}%
satisfying, $Tr\circ \overline{\omega }_{\alpha }=\omega _{\alpha }$.

Now it is clear that all the decompositions (\ref{decom1}) ,(\ref{dec2}),( %
\ref{decom3}), and (\ref{decom4}) are compatible with the action of the
Galois group $G$, hence, they induce corresponding decompositions on the
level of invariants, which establish the content of the lemma.

We are left to show that for each $\alpha \in \Xi $, the algebra $K_{\alpha
} $ is a field, and, in addition, $\dim _{K_{\alpha }}V_{\alpha }=2$. The
second claim follows from the (very easy to verify) fact that the vector
space $\mathbf{V}_{\alpha }$ is a free module of rank $2$ over the algebra $%
\mathbf{A}_{\alpha }$. In proving the first claim, we individually analyze
two cases

\begin{itemize}
\item The orbit $\mathcal{O}_{\alpha }$ is of type I. In this case $%
T_{\alpha }$ acts irreducibly on $V_{\alpha }$ and we are back in the
completely inert situation of Subsection \ref{particular}.

\item The orbit $\mathcal{O}_{\alpha }$ is of type II. \ In this case $%
V_{\alpha }=V_{\beta }\oplus V_{\widehat{\beta }}$ where $V_{\beta },V_{%
\widehat{\beta }}$ are irreducible $T_{\alpha }$ invariant Lagrangian
subspaces. Respectively, the algebra $A_{\alpha }$ decomposes into a direct
sum of two fields $A_{\alpha }=A_{\beta }\oplus A_{\widehat{\beta }}$. Note
that $A_{\beta }=Z(T_{\alpha },End(V_{\beta }))$ and $A_{\widehat{\beta }%
}=Z(T_{\alpha },End(V_{\widehat{\beta }}))$, which implies that $A_{\beta
},A_{\widehat{\beta }}$ are indeed fields. Moreover, the involution $\Theta $
induces an isomorphism $A_{\beta }\simeq A_{\widehat{\beta }}$ and $%
K_{\alpha }$ is the diagonal field with respect to this isomorphism.
\end{itemize}

This concludes the proof of Lemma \ref{technical}.

\subsection{Proof of Proposition \protect \ref{resHeisenberg}\label%
{PresHeisenberg}}

The proof is immediate. Clearly, the representation $\overline{\pi }$ is
irreducible (the homomorphism $\iota _{H}$ is surjective). Therefore, by the
Stone--von Neumann theorem (Theorem \ref{SVN}) it is determined by its
central character. The homomorphism $\iota $ restricts \ to an homomorphism $%
K=Z(\overline{H})\rightarrow Z(H)=k$ between the centers, which is given by $%
Tr=Tr_{K/k}$. Consequently, this implies that $\overline{\pi }_{|Z(\overline{%
H})}=$ $\overline{\psi }\cdot Id$. Concluding the proof of the proposition.

\subsection{Proof of Theorem \protect \ref{reWeil}\label{PreWeil}}

It is enough to show that the representations $\overline{\rho }$ and $%
\overline{\pi }$ satisfy the compatibility condition (\ref{Egorov1}).
Indeed, for $g\in \overline{Sp}$ and $h\in \overline{H}$ 
\begin{eqnarray*}
\overline{\rho }(g)\overline{\pi }(h)\overline{\rho }(g)^{-1} &=&\rho (\iota
_{S}(g))\pi (\iota _{H}(h))\rho (\iota _{S}(g))^{-1} \\
&=&\pi (\iota _{S}(g)\cdot \iota _{H}(h))=\pi (\iota _{H}(g\cdot h)=%
\overline{\pi }(g\cdot h),
\end{eqnarray*}%
concluding the proof.

\subsection{Proof of Proposition \protect \ref{Heisenbergdecom}\label%
{PHeisenbergdecom}}

The representation $\overline{\pi }$ is irreducible ($\iota _{H\text{ }}$is
surjective), this immediately implies that $(\overline{\pi },\overline{H},%
\mathcal{H)\simeq (}\otimes \overline{\pi }_{\alpha },\Pi \overline{H}%
_{\alpha },\otimes \mathcal{H}_{\alpha })$, where $\overline{\pi }_{\alpha }$
is irreducible, for every $\alpha $. By the Stone--von Neumann theorem a
representation of the Heisenberg group $\overline{H}_{\alpha }$ is
characterized by the action of the center, which clearly in this case acts
via the character $\overline{\psi }_{\alpha }=\psi \circ Tr_{K_{\alpha }/k}$%
. This completes the proof.

\subsection{Proof of Theorem \protect \ref{Weildecomp}\label{PWeildecomp}}

The representation $\overline{\rho }$ is characterized \ by the identity (%
\ref{Egorov1}) with respect to the representation $\overline{\pi }$. Using
Proposition \ref{Heisenbergdecom}, it is enough to show that $\otimes 
\overline{\rho }_{\alpha }$ satisfies the identity (\ref{Egorov1}) with
respect to the representation $\otimes \overline{\pi }_{\alpha }$, which is
principally reduced to the case analyzed in Subsection \ref{PreWeil}. This
completes the proof.

\subsection{Proof of Proposition \protect \ref{CT}\label{PCT}}

Given the realization of the symplectic rank $r_{p}=\left \vert \Xi
\right
\vert $ that appears in Subsection \ref{Ptechnical}, the proof of
Proposition \ref{CT} is an immediate consequence of the Chebotarev density
theorem \cite{Se}.

\section{Proof of Theorem \protect \ref{twodim1}\label{Pbound}}

We need to show that $\left \vert \sum_{g\in T}\chi ^{-1}(g)ch_{\tau
}(g,v)\right \vert \leq 2\sqrt{q}.$ It is easy to verify that for $g=1,$ $%
v\neq 0$ we have $ch_{\tau }(g,v)=0.$ Moreover, for $g\neq 1$ using the
character formula (\ref{Char H-W}) and Corollary \ref{equ} we have $ch_{\tau
}(g,v)=\pm \sigma _{T}(g)\psi (\tfrac{1}{2}\omega (\frac{v}{g-I},v))$ where $%
\sigma _{T}$ is the unique quadratic character of $T.$ We may assume $\chi
\neq $ $\sigma _{T}$ (the case $\chi =$ $\sigma _{T}$ being trivial).
Concluding, \ for a non-trivial character $\chi $ of $T$ we consider the
function $f_{\chi }(g)=\chi (g)\psi (\tfrac{1}{2}\omega (\frac{v}{g-I},v))$
on the set $X=T\smallsetminus I$ and we would like to show that 
\begin{equation}
\left \vert c_{\chi }=\sum_{g\in X}f_{\chi }(g)\right \vert \leq 2\sqrt{q}.
\label{Est}
\end{equation}

\subsection{Solution via geometrization}

Our problem fits nicely into Grothendieck's geometrization methodology:
replacing sets by varieties; functions be sheaves; and reducing the proof of
the estimate (\ref{Est}) to a topological statement on certain cohomology
groups.

\subsubsection{Replacing sets by varieties}

We denote by $\overline{k}$ an algebraic closure of the field $k=\mathbb{F}%
_{q}$. We use bold-face letters to denote a variety $\mathbf{Y}$ and normal
letters $Y$ to denote its corresponding set of rational points $Y=\mathbf{Y}%
(k)$. By a variety $\mathbf{Y}$ over $k$ we mean a quasi-projective
algebraic variety, such that the defining equations are given by homogeneous
polynomials with coefficients in the finite field $k$. In this situation,
there exists a (geometric) \textit{Frobenius} endomorphism $Fr:\mathbf{%
Y\rightarrow Y}$, which is a morphism of algebraic varieties. We denote by $%
Y $ the set of points fixed by $Fr$, i.e., $Y=\mathbf{Y}(k)=\{y\in \mathbf{Y}%
:Fr(y)=y\}.$ \ 

Important example for us are $k$ and $k^{\ast }$ which are the sets of
rational points of the algebraic varieties $\mathbb{G}_{a}$ and $\mathbb{G}%
_{m},$ respectively, and $T$ and $X=T\smallsetminus I$ which are the sets of
rational points of the algebraic varieties $\mathbf{T\subset Sp}$ and $%
\mathbf{X=T\smallsetminus I.}$

\subsubsection{Replacing functions by sheaves}

Let $\mathsf{D}^{b}(\mathbf{Y)}$ denote the bounded derived category of
constructible $\ell $-adic sheaves on $\mathbf{Y}$ \cite{BBD, D1, KW}. A
Weil structure associated to an object $\mathcal{F\in }\mathsf{D}^{b}(%
\mathbf{Y)}$ is an isomorphism $\theta :Fr^{\ast }\mathcal{F}\overset{\sim }{%
\longrightarrow }\mathcal{F}$. A pair $(\mathcal{F},\theta )$ is called a
Weil object. By an abuse of notation we often denote $\theta $ also by $Fr$.
We choose once an identification $\overline{%
%TCIMACRO{\U{211a} }%
%BeginExpansion
\mathbb{Q}
%EndExpansion
}_{\ell }\simeq 
%TCIMACRO{\U{2102} }%
%BeginExpansion
\mathbb{C}
%EndExpansion
$, hence all sheaves are considered over the complex numbers. Given a Weil
object $(\mathcal{F},Fr^{\ast }\mathcal{F\simeq F})$ one can associate to it
a function $f^{\mathcal{F}}:Y\rightarrow 
%TCIMACRO{\U{2102} }%
%BeginExpansion
\mathbb{C}
%EndExpansion
$ to $\mathcal{F}$ as follows 
\begin{equation}
f^{\mathcal{F}}(y)=\tsum \limits_{i}(-1)^{i}Tr(Fr_{|H^{i}(\mathcal{F}_{y})}).
\label{S-t-F}
\end{equation}%
This procedure is called \textit{Grothendieck's sheaf-to-function
correspondence \cite{Gr}}.

Important examples for us are the additive character $\psi :k\longrightarrow 
%TCIMACRO{\U{2102} }%
%BeginExpansion
\mathbb{C}
%EndExpansion
^{\times }$ which is associated via the sheaf-to-function correspondence to
the Artin--Schreier sheaf $\mathcal{L}_{\psi }$ on the variety $\mathbb{G}%
_{a}$, i.e., we have $f^{\mathcal{L}_{\psi }}=\psi ;$ the Legendre character 
$\sigma $ on $k^{\ast }$ which is associated to the Kummer sheaf $\mathcal{L}%
_{\chi }$ on the variety $\mathbb{G}_{m}$, \ i.e., $f^{\mathcal{L}_{\chi
}}=\chi $; and the function $f_{\chi }$ on $X$ which is associated with the
Weil sheaf $\mathcal{F}_{\chi }\mathcal{=L}_{\chi (g)}\otimes \mathcal{L}%
_{\psi (1/2\omega (\frac{v}{g-1},v))}$ on $\mathbf{X.}$ Finally, we have the
Weil object%
\begin{equation*}
\mathcal{C}_{\chi }=\tint \mathcal{F}_{\chi }\in \mathsf{D}^{b}(\mathbf{pt}),
\end{equation*}%
where $\int =\int_{!}$ denotes integration with compact support \cite{D1}.
The Grothendieck--Lefschetz trace formula \cite{Gr} implies that $f^{%
\mathcal{C}_{\chi }}=c_{\chi }.$

\subsubsection{Geometric statement}

The sheaf $\mathcal{F}_{\chi }$\ is a non-trivial rank $1$\ irreducible
local system of pure weight zero $w(\mathcal{F}_{\chi })=0$. The deep
theorem of Deligne \cite{D1} say that integration with compact support does
not increase weight; therefore $\mathcal{C}_{\chi }$ is of mixed weight $w(%
\mathcal{C}_{\chi })\leq 0,$ i.e., we have the following bound on each
eigenvalue $\lambda $ of $Fr:$ 
\begin{equation*}
|\lambda (Fr_{|H^{i}(\mathcal{C}_{\chi })})|\leq \sqrt{q}^{i}.
\end{equation*}%
This means by (\ref{S-t-F}) that in order to get the bound (\ref{Est}) it is
enough to prove the following geometric statement:

\begin{lemma}[Vanishing lemma]
The object $\mathcal{C}_{\chi }$ is cohomologically supported at degree $1$
and in addition $\mathrm{dim}$ $H^{1}($ $\mathcal{C}_{\chi })=2$.
\end{lemma}

\subsubsection{Proof of the Vanishing lemma}

The fact that only $H^{1}($ $\mathcal{C}_{\chi })$ does not vanishes follows
from the fact that the local system $\mathcal{F}$ is of rank $1$ and
non-trivial, i.e., it admits non-trivial monodromy. We are left to compute
the dimension of $H^{1}(\mathcal{C}_{\chi }).$ Because $\mathcal{C}_{\chi }$
is cohomologically supported at degree $1$ we have 
\begin{equation*}
\mathrm{dim}H^{1}(\mathcal{C}_{\chi })=-\chi _{Fr}(\mathcal{C}_{\chi }%
\mathcal{)},
\end{equation*}%
where $\chi _{Fr}$ denotes the Euler characteristic. This is a topological
invariant defined by $\chi _{Fr}(\mathcal{C}_{\chi })=\tsum
\limits_{i}(-1)^{i}\dim $ $H^{i}(\mathcal{C}_{\chi }).$

The actual computation of the Euler characteristic $\chi _{Fr}(\mathcal{C}%
_{\chi }\mathcal{)}$ is done using the \textit{%
Ogg--Shafarevich--Grothendieck formula} \cite{K, KW}. Recall that $\mathcal{C%
}_{\chi }=\tint \mathcal{F}_{\chi }.$ We have 
\begin{equation}
\chi _{Fr}(\tint \overline{{\mathbb{Q}}}_{\ell })-\chi _{Fr}(\tint \mathcal{F%
}_{\chi })=\sum_{y\in \mathbf{Y}\setminus \mathbf{X}}\mathrm{Swan}_{y}(%
\mathcal{F}_{\chi }),  \label{OSG}
\end{equation}%
where $\overline{{\mathbb{Q}}}_{\ell }$ denotes the constant sheaf on $%
\mathbf{X}$ and $\mathbf{Y}$ is some compact curve containing $\mathbf{X}$.
This formula expresses the difference of $\chi _{Fr}(\tint \limits_{\mathbf{X%
}}\overline{{\mathbb{Q}}}_{\ell })$ from $\chi _{Fr}(\tint \limits_{\mathbf{X%
}}\mathcal{F})$ as a sum of local contributions. We take $\mathbf{Y}=\mathbb{%
P}^{1};$ therefore under this choice the complement $\mathbf{Y\backslash X}$
consists of 3 points\ $\mathbf{Y\backslash X=}\left \{ 0,1,\infty \right \} $%
. \ Using the standard properties of the Swan conductors (see \cite{K, KW}),
and the well known values of the Swan conductors of the standard sheaves $%
\mathcal{L}_{\chi \text{ }}$and $\mathcal{L}_{\psi }$, we obtain $Swan_{0}(%
\mathcal{F)}=Swan_{0}(\mathcal{L}_{\chi \text{ }})=0,$ $Swan_{1}(\mathcal{F)}%
=Swan_{\infty }(\mathcal{L}_{\psi \text{ }})=1,$ and $Swan_{\infty }(%
\mathcal{F)}=Swan_{\infty }(\mathcal{L}_{\chi \text{ }})=0.$ In addition, we
have that $\chi _{Fr}(\tint \overline{{\mathbb{Q}}}_{\ell })=-1$.
Concluding, using formula (\ref{OSG}) we get $\chi _{Fr}(\tint \limits_{%
\mathbf{X}}\mathcal{F})=-2.$

This completes the proof of the Vanishing lemma and the proof of Theorem \ref%
{twodim1}.

\bigskip

\end{document}